\documentclass[a4,11pt]{article}
\usepackage{graphicx,graphics,epsfig}
\usepackage{amsmath,amsthm,amssymb}%,dsfont}
\usepackage[usenames]{color}
\usepackage{subfig}
\graphicspath{{./pics/}}
\newtheorem{Prop}{Proposition}
\newtheorem{lem}{Lemma}
\newtheorem{thm}{Theorem}
\newtheorem{cor}{Corollary}%[thm]

\newtheorem{rem}{Remark}

\allowdisplaybreaks[2]
\newcommand{\R}{\mathbb{R}}
\newcommand{\N}{\mathbb{N}}
\newcommand{\Z}{\mathbb{Z}}

\title{CIP methods for hyperbolic system with variable and discontinuous coefficient}

\author{Kazufumi Ito
\thanks{Center for Research in Scientific Computation \& Department of Mathematics,
North Carolina State University,
Raleigh,
NC 27695, USA.
(kito@math.ncsu.edu)}
\and
Tomoya Takeuchi
\thanks{Center for Research in Scientific Computation,
North Carolina State University,
Raleigh,
NC 27695, USA.
({\tt tntakeuc@ncsu.edu})}
}

\begin{document}
\maketitle
\begin{abstract}
We propose a multi-moment method for one-dimensional hyperbolic
equations with smooth coefficient and piecewise constant coefficient.
The method is entirely based on the backward characteristic method and uses the solution and its
derivative as unknowns and cubic Hermite interpolation for each
computational cell.
The exact update formula for solution and its
derivative is derived and used for an efficient time integration.
At points of discontinuity of wave speed we define a
piecewise cubic Hermite interpolation based on immersed interface method.
The method is extended to the one-dimensional Maxwell's equations with
variable material properties.
\end{abstract}
%\begin{keywords}
%CIP method and exact time integration, hyperbolic equations with smooth coefficient, hyperbolic equations in discontinuous media.
%\end{keywords}

%\begin{AMS}
%\end{AMS}
%
%\pagestyle{myheadings}
%\thispagestyle{plain}
%\markboth{K. ITO AND T. TAKEUCHI}{CIP FOR VARIABLE COEFFICIENT}
%\newtheorem{rem}{Remark}
\section{Introduction}
In this paper we develop a multi-moment method for the wave
propagation, for instance,
\begin{equation} \label{pde}
u_t + (c(x) u)_x = 0, \quad t>0, \ x\in \R,\qquad u(0,x)=u_0(x),\;\;
x\in \R.
\end{equation}

Multi-moment methods approximate the solutions to differential equations using not only the
primitive variable but also another numerical information at the grid or the cell such as its derivatives, the cell average of the numerical solutions. A Hermite polynomial is usually defined to interpolate such information. Then the numerical quantities are evolved in time simultaneously.

Our method is strongly motivated by and
closely related to CIP methods, which is one of the multi-moment methods.
It was first proposed in \cite{Takewaki+NishiguchiETAL-Cubiintepseumeth:85} for constant velocity field. The CIP method uses the exact integration in time by the characteristic method and
uses the cubic Hermite polynomial in each cell $[x_{j-1},x_j]$ based on solution values and its derivatives at two endpoints
$x_{j-1},\;x_j$. The method provides an accurate, less-dispersive and less-dissipative numerical solution.
Here is an incomplete list of references for CIP method and related works:
nonlinear hyperbolic equations \cite{Yabe+Aoki-univsolvhypeequa:91},
multi dimensional hyperbolic equations \cite{Yabe-univcubiintesolv:91,Yabe+IshikawaETAL-univsolvhypeequa:91},
the multi-phase analysis \cite{Yabe+XiaoETAL-consinteprofmeth:01},
a multi-dimensional the Maxwell's equations \cite{Ogata+YabeETAL-AccuNumeScheMaxw:06}, light propagation in dielectric media \cite{Barada+FukudaETAL-Cubiintepropsche:06},
a new mesh system applicable to non-orthogonal coordinate system \cite{Yabe+MizoeETAL-Highschewithmeth:04},
a numerical investigation of the stability and the accuracy \cite{Utsumi+KunugiETAL-StabaccuCubiInte:97}.
The other method closely related to CIP, we refer to \cite{Aoki-Intediffoper:97,Imai+Aoki-highimplscheappl:06}.

Our contributions are as follows:
Firstly, we develop the exact solution formula for solution and its derivative for \eqref{pde}
with smooth variable wave speed $c(x)$.
The cubic Hermite polynomial is then used to evaluate the formula approximately.
Our development of the CIP scheme is entirely based on the characteristic
method and results in a different (improved and simpler) CIP scheme than
the conventional one for equations with variable wave speed. For
example the new scheme allows us to take an arbitrary time step (no
CFL limitation) without losing the stability and accuracy.
%We investigate the convergence of the proposed method for a constant speed.

Next, we develop an immersed interface method \cite{Li+Ito-immeintemeth:06,Zhang+LeVeque-immeintemethacou:97} for \eqref{pde} with a piecewise constant wave speed. We define a piecewise cubic Hermite polynomial on a cell $[x_{j-1},x_j]$ which contains a point of jump discontinuity in $c(x)$ using proper interface conditions and solution values and its derivatives at two endpoints $x_{j-1},\;x_j$.

Lastly, this interface treatment is applied to the one-dimensional Maxwell's equation with variable material properties:
We first approximate the variable coefficients by a piecewise constant (discontinuous) coefficient.
The d'Alembert's based method for the Maxwell's equations is developed for the piecewise constant media and then applied to Maxwell system with piecewise constant coefficients.

An outline of our presentation is: in Section \ref{sec:CIPmethod} our proposed method for hyperbolic equations with smooth variable coefficient is developed, in Section \ref{sec:error} the error analysis of the CIP scheme for the constant velocity is presented, in Section \ref{sec:IIMCIP} a CIP scheme for discontinuous media is developed, and in Section \ref{sec:IIMCIPmax} its application to the Maxwell's system is presented.
Each section contains some numerical tests to verify the accuracy of the proposed method.
\section{CIP method for smooth coefficient}\label{sec:CIPmethod}
In this section we propose a CIP method for hyperbolic equations with smooth coefficient.
%We first develop the exact time integration for the solutions of hyperbolic equations.
We consider the advection equation \eqref{pde} with sufficiently smooth $c(x)$ as a model equation. The characteristic method is the key for designing a CIP scheme.
%CIP (Cubic Interpolated Pseudo-particle) \\
%CIP (Cubic Interpolated Propagation) \\
%CIP (Constrained Interpolation Profile) \\
%CIP(Constrained Interpolated Profile)\\
%\subsection{CIP overview}
Let $u$ be the exact solution of the equation. We discretize the
time domain and the spatial domain with grid size $\Delta t>0 $ and
$\Delta x>0$. Write $t_n=n\Delta t$ and $x_k = k \Delta x$, $n\in \N$, $ k\in
\Z$. Let us consider the characteristic curve $x=x(t)$ subject to $x(\Delta t)=x_k$ to the equation \eqref{pde},
\begin{equation}
  \frac{d x(t)}{dt} = c(x(t)),\quad x(\Delta t)=x_k. \label{ode:x}
\end{equation}
We integrate the equation backward in time to find $y_k=x(0)$. The following update formula is the fundamental for developing CIP scheme.
\begin{Prop}
The exact solution and its derivative at time $t=t_{n+1}$, $x=x_k$ are
\begin{equation}\label{MOC}
\begin{array}{l}
u(x_k,t+\Delta t)=\dfrac{c(y_k)}{c(x_k)}u(y_k,t), \\
u_x(x_k, t+\Delta t) =
\left(\dfrac{c'(y_k)}{c(x_k)}-\dfrac{c'(x_k)}{c(x_k)}\right)
\dfrac{c(y_k)}{c(x_k)} u(y_k,t) + \left(\dfrac{c(y_k)}{c(x_k)} \right)^2 u_x(y_k,t).
\end{array}
\end{equation}
\end{Prop}
\begin{proof}
%The system of characteristic equations for \eqref{pde} is
%comprised of (see, e.g., p 98 in \cite{Evans-Partdiffequa:98})
Let $p_1(t)=u_t(t_n+t,x(t))$, $p_2(t)=u_x(t_n+t,x(t))$ and $z(t) =
u(t_n+t,x(t))$. The system of ODE for $p_1,\;p_2, \;z$ becomes (\cite{Evans-Partdiffequa:98}, p. 98)
\begin{align}
& p_1(t) + c(x(t))p_2(t) + c^\prime(x(t))z(t) = 0, \label{1storder} \\
& \frac{d p_1(t)}{dt} = -c^\prime(x(t))p_1(t), \label{ode:p1}\\
% \frac{d p_2(t)}{dt} &=c_x(x(t))p_2(t)+c_{xx}(x(t))-c^\prime(x(t)), \label{ode:p2}\\
&  \frac{d z(t)}{dt} = p_1(t)+c(x(t))p_2(t). \label{ode:z}
%&  \frac{d x(t)}{dt} = c(x(t)), \label{ode:x}
\end{align}
%Firstly, we integrate \eqref{ode:x} subject to $x(\Delta t)=x_k$ backward in
%time to find $y_k=x(0)$. \\
From \eqref{1storder} and \eqref{ode:z}, $\frac{d z(t)}{dt}= -c^\prime(x(t)) z(t).$ Hence, from \eqref{ode:x}
\begin{equation*}
\frac{d}{dt}(c(x(t)) z(t)) = c^\prime(x(t))\frac{d}{dt}x(t)z(t) +
c(x(t)) \frac{d}{dt}z(t) = c(x(t))\left(z(t)c^\prime(x(t))+\frac{d
}{dt}z(t)\right) = 0.
\end{equation*}
Upon integrating this, we obtain $z(t)=\frac{c(y_k)}{c(x(t))}z(y_k).$ Similarly, from \eqref{ode:p1}
%\begin{equation*}
%  \int_{0}^{s} \frac{\frac{d }{d\tau}p_1(\tau)}{p_1(\tau)}\,d\tau
%  = -\int_0^s c_x(x(\tau))\,d\tau
%  = -\int_{x(s)}^{x(0)} \frac{c_x(x(\tau))}{\frac{d}{d\tau}x(\tau) }\,d\tau
%\end{equation*}
%Hence
$
p_1(t)=\frac{c(y_k)}{ c(x(t))}p_1(y_k),
$ and from \eqref{1storder}
\begin{equation*}
\dfrac{c(y_k)}{c(x(t))}p_1(y_k)+c(x(t))p_2(x(t)+c^\prime(x(t))\dfrac{c(y_k)}{c(x(t))}z(y_k)=0
\end{equation*}
where
$p_1(y_k)=-c(y_k)p_2(y_k)-c^\prime(y_k)z(y_k).$
Hence, we obtain
\begin{equation*}
p_2(t) =
\left(\dfrac{c'(y_k)}{c(x(t))}-\dfrac{c'(x(t))}{c(x(t))}\right)
 \dfrac{c(y_k)}{c(x(t))}z(y_k) + \left(\dfrac{c(y_k)}{c(x(t))} \right)^2 p_2(y_k).
\end{equation*}
Thus, we have \eqref{MOC}.
\end{proof}

The formula \eqref{MOC} is the fundamental relations and an explicit numerical scheme is constructed based on the relation. Both of the primitive variable
$u$ and the spatial derivative $u_x$ will be simultaneously updated in CIP scheme.

With \eqref{MOC} we can compute the values
$u(t_{n+1},x_k)$ and $u_x(t_{n+1},x_k)$ if we know the values
$u(t_n,y_k)$ and $u_x(t_n,y_k)$. Suppose  $y_k$ falls in the $j=j(k)$-th interval $(x_{j-1},x_{j})$.
%We interpolate the value $u(t_n,y_k)$ in terms of known numerical quantities at grid points $x_{j-1},\;x_j$
%at time level $t=t_n$.
Let us denote the numerical solution for $u(t_n,x_j)$ and $u_x(t_n,x_j)$ by $u^n_j$ and $v^n_j$ respectively, and suppose that they are given at all nodes as variables. We construct
the cubic Hermite polynomial on the interval
$[x_{j-1},x_j]$;
\begin{equation}\label{cubic}
\begin{array}{l}
H_j(x)=
u^n_{j-1}p_1( \frac{x-x_{j-1}}{\Delta x}) + u^n_{j}p_2( \frac{x-x_{j-1}}{\Delta x}) \\
\qquad + \Delta x\; v^n_{j-1}q_1( \frac{x-x_{j-1}}{\Delta x}) + \Delta x\;v^n_{j}q_2( \frac{x-x_{j-1}}{\Delta x}),
     \end{array}
\end{equation}
where $p_1(\xi) = (\xi-1)^2(2\xi+1),$ $p_2(\xi) = \xi^2(3-2\xi),$
 $ q_1(\xi) = (\xi-1)^2\xi,$ $q_2(\xi) = \xi^2(\xi-1).$
Approximating  $u(t_n,y_k)\approx H_j(y_k)$ and $v(t_n,y_k)\approx
H_j^\prime(y_k)$ on the interval $(x_{j-1},x_j)$ and from
\eqref{MOC}, we arrive at CIP scheme:
\begin{equation}\label{newu}
  \begin{array}{ll}
  u^{n+1}_k&= \frac{c(y_k)}{c(x_k)}H_j(y_k),\\
  v^{n+1}_k &=
  \frac{c(y_k)}{c(x_k)}\left[\frac{c'(y_k)}{c(x_k)}-\frac{c'(x_k)}{c(x_k)}\right]H_j(y_k)
  +\left(\frac{c(y_k)}{c(x_k)}\right)^2 H_j^\prime(y_k).
  \end{array}
\end{equation}
The cubic polynomial $H_j$ associated with the numerical solution is called \emph{cubic interpolation
profile} on the interval $[x_{j-1},x_{j}]$.

For the transport equation $u_t+c(x)u_x=0$, we have the (exact) update formula:
\begin{equation}\label{tra}
  u^{n+1}_k= H_j(y_k),\quad
  v^{n+1}_k =
  \frac{c(y_k)}{c(x_k)}H_j^\prime(y_k).
\end{equation}

 In summary, CIP scheme
for \eqref{pde} is composed of
\begin{description}
  \item[CIP0] Solve the characteristic ODE \eqref{ode:x} subject to $x(\Delta t)=x_k$ to find $y_k=x(0)$.
  \item[CIP1] For each time level $t_n$ construct the cubic interpolation profile $H_j(x)$.
  \item[CIP2] Update $u^{n+1}_k$, $v^{n+1}_k$ by \eqref{newu}.
\end{description}

Suppose that the velocity filed $c$ is constant, there exist an integer $\ell$ such that $y_k = x_{k}- c\Delta t \in [x_{k-1-\ell},x_{k-\ell}]$ for all $k$ and thus the quantity $\lambda:=\frac{x_{k-\ell}-y_k}{\Delta x}$ is independent of $k$ and $0\le \lambda \le 1$, and the one step map \eqref{newu} becomes, denoting $k-\ell$ by $j$,
\begin{equation}\label{u}
\begin{array}{ll}
u^{n+1}_k=
%u^n_{j-1}p_1(1-\lambda) + u^n_{j}p_2(1-\lambda) +  \Delta x\; v^n_{j-1}q_1(1-\lambda) + \Delta x\; v^n_{j}q_2(1-\lambda)
%\\ \\
&\lambda^2 ( 3-2 \lambda) u^n_{j-1} + (1-\lambda)^2(1+2 \lambda) u^n_{j}\\ \\
&+\Delta x\;\lambda ^2(1-\lambda)v^n_{j-1}  + \Delta x (1-\lambda)^2(-\lambda)\; v^n_{j},\\
\\
v^{n+1}_k=
%u^n_{j-1}p'_1(1-\lambda) + u^n_{j}p'_2(1-\lambda) +  \Delta x\; v^n_{j-1}q'_1(1-\lambda) + \Delta x\; v^n_{j}q'_2(1-\lambda)
&\frac{-6 \lambda(1-\lambda)}{\Delta x} u^n_{j-1} +  \frac{6 \lambda(1-\lambda)}{\Delta x} u^n_{j}\\ \\
& + ( 3\lambda -2)\lambda v^n_{j-1} + (1-\lambda)(1-3\lambda) v^n_{j} .
  \end{array}
\end{equation}
We shall investigate the consistency and the stability of the scheme in the next section.

In existing articles on CIP, the approximate scheme \eqref{newu} is used in the following manner. One differentiates
with respect to $x$ the non-conservative form $ u_t + cu_x = -c' u $
of the equation \eqref{pde} to obtain $ v_t + cv_x =(-c' u)_x- c'
u_x$, here we have used the notation $v=u_x$. Next, these equations
are split into the advection phase $u_t + cu_x=0 ,\quad v_t +
cv_x=0$, and the non-advection phase $ u_t = - c' u,\quad   v_t
=(-c' u)_x- c' u_x$. The numerical solution of the advection phase
are given using the cubic Hermite polynomial $u^*_k=H(x^*)$,
$v^*_k=H_x(x^*)$, where the upwind point $x^\ast$ is determined appropriately for each $k$.
Then the system of the non-advection equation is
integrated to obtain the numerical solution of the next time level.
In \cite{Yabe+Aoki-univsolvhypeequa:91}, the following approximate scheme for
the equations is proposed:
\begin{equation}\label{sol1}
\begin{array}{ll}
 u^{n+1}_k &= (1-c'(x_k)\Delta t)u^*_k,  \\
 v^{n+1}_k &= \left(1-\tfrac{c(x_{k+1})-c(x_{k-1})}{2\Delta x}\right)v^*_k +
  \tfrac{u^*_{k+1}-u^*_{k-1} -u^n_{k+1}+u^n_{k-1} }{2\Delta x}.
 \end{array}
\end{equation}
The other variant is proposed in \cite{Yabe+MizoeETAL-Highschewithmeth:04}:
\begin{equation}\label{sol2}
 u^{n+1}_k = \exp(-c'(x_k)\Delta t)u^*_k,  \quad
 v^{n+1}_k = -u^*_kc''(x_k)\Delta t + \exp(-c'(x_k)\Delta t)v^*_k.
\end{equation}
Note that the CFL number is limited less than or equal to 1 in these methods.
%and the splitting procedure causes error in time and the method reduces to second order accuracy at most.
One can verify that \eqref{sol1} and \eqref{sol2} are approximations to \eqref{newu}. Indeed, if $x_k - y_k\le \Delta x$ is sufficiently small, we have $\tfrac{c(y_k)}{c(x_k)}\approx \tfrac{c(y_k-c(x_k)\Delta
t)}{c(x_k)}\approx 1-c'(y_k)\Delta t \approx \exp(-c'(x_k)\Delta t)$.

\subsection{Error analysis}\label{sec:error}
In this section we investigate the consistency and the stability of the CIP scheme \eqref{u} for constant velocity and derive an error estimates for the numerical solution. Utsumi et al \cite{Utsumi+KunugiETAL-StabaccuCubiInte:97} studies the stability and accuracy of the scheme through numerical simulation.
We give a theoretical underpinning for the numerically observed evidence in terms of the eigenvalue analysis.

The one step map \eqref{u} is written
\begin{equation*}
\left(\begin{array}{c} u^{n+1}_k \\ \tilde{v}^{n+1}_k
\end{array}\right) =A(\lambda)\left(\begin{array}{c} u^{n}_j \\
\Delta x\;{v}^{n}_j
\end{array}\right)+B(\lambda)\left(\begin{array}{c} u^{n}_{j-1} \\ \Delta x\;{v}^{n}_{j-1} \end{array}\right)
\end{equation*}
where $\lambda=\dfrac{x_j-y_k}{\Delta x}$ ($0 \le\lambda\le 1$), and
\begin{equation*}
A(\lambda)=\left(\begin{array}{cc} (1-\lambda)^2(1+ 2\lambda) & -
(1-\lambda)^2\lambda \\ \\
6\lambda(1-\lambda) & (1-\lambda)(1-3 \lambda)
\end{array}\right),\quad
B(\lambda)=\left(\begin{array}{cc} \lambda ^2(3-2 \lambda) &
\lambda ^2(1- \lambda) \\ \\ -6 \lambda(1-\lambda) & (3\lambda-2)\lambda
\end{array}\right).
\end{equation*}
We recall that there exist an integer $\ell$ such that $j=k-\ell$ for all $k$ since the velocity is constant.
%\begin{Prop}
%Suppose $\sum_{k} v^0_k=0$, then $\sum_k u^n_k=\sum_k u^0_k$ for
%all $n\ge0$. \end{Prop}
%\begin{proof} It follows from
%$$
%A(\lambda)+B(\lambda)=\left(\begin{array}{cc} 1 &
%-\lambda+3\lambda^2 -2\lambda ^3 \\  0 & 1-6\lambda +6  \lambda ^2
%\end{array}\right).
%$$
%\end{proof}

Let us denote the exact solution of \eqref{pde} for constant coefficient at the node $x_k$ at time $t_n$ by $U^n_k=u(t_n,x_k)$ and $V^n_k=u_x(t_n,x_k)$.
\begin{thm}
\begin{equation*}
\left(\begin{array}{c} U^{n+1}_k \\ \Delta x\; V^{n+1}_k
\end{array}\right)= A(\lambda)\left(\begin{array}{c} U^{n}_j \\
\Delta x\;{V}^{n}_j
\end{array}\right)+B(\lambda)\left(\begin{array}{c}  U^{n}_{j-1} \\ \Delta x \;{V}^{n}_{j-1} \end{array}\right)
+O((\Delta x)^4).
\end{equation*}
Therefore, CIP scheme \eqref{u} is consistent and is of order four in space for $u(x,t)$ and of order three for $u_x(x,t)$.
 \end{thm}
 \begin{proof}
This estimate directly follows from the point wise estimate for Hermite interpolation. Indeed from \cite{Agarwal-SharHerminteerro:91}, denoting the Hermite polynomial of $u(t_n,x)$ on $[x_{j-1},x_j]$ by $h(x)$, we have
\begin{equation*}
\begin{array}{ll}
  u(t_n,x)-h(x)  &= \frac{u^{(4)}(t_n,\xi_0)}{4!}(x-x_{j-1})^2(x-x_j)^2, \\
  u_x(t_n,x)-h_x(x) &= \frac{u^{(4)}(t_n,\xi_1)}{3!}(x-x_{j-1})(x-x_j)(x-\xi_2),
\end{array}
\end{equation*}
for $x\in [x_{j-1},x_j]$, where $\xi_0,\;\xi_1, \;\xi_2\in (x_{j-1},x_j)$.
 Thus, we have
\begin{align*}
\left(\begin{array}{c}
U^{n+1}_k \\ \Delta x\; V^{n+1}_k
\end{array}\right)&=
\left(\begin{array}{c}
u(t_n,y_k) \\ \Delta x\; u_x(t_n,y_k)
\end{array}\right)
=
\left(\begin{array}{c}
 h(y_k)+ (E^{1})^n_j  \\ \Delta x\;(h_x(y_k)+ (E^{2})^n_j)
\end{array}\right) \\
& =
A(\lambda)
\left(\begin{array}{c} U^{n}_j \\
\Delta x\;{V}^{n}_j
\end{array}\right)
+B(\lambda)\left(\begin{array}{c}  U^{n}_{j-1} \\ \Delta x \;{V}^{n}_{j-1} \end{array}\right)
+\left(\begin{array}{c}
(E^{1})^n_j  \\ \Delta x\; (E^{2})^n_j
\end{array}\right),
\end{align*}
where $j=k-\ell$ and
\begin{equation*}%\label{E1E2}
\begin{array}{lll}
(E^{1})^n_j &=\tfrac{u^{(4)}(t_n,\xi_{0})}{4!}\lambda^2(1-\lambda)^2\,(\Delta
x)^4 = \tfrac{u^{(4)}(0,\xi_{0} - c n\Delta t)}{4!}\lambda^2(1-\lambda)^2\,(\Delta x)^4\ \\
(E^{2})^n_j &= \tfrac{u^{(4)}(t_n,\xi_{1})}{3!} \lambda(1-\lambda)(\lambda+\tfrac{\xi_{2}}{\Delta x})(\Delta x)^3 \\ &=
\tfrac{u^{(4)}(0,\xi_{1}- c n\Delta t)}{3!} \lambda(1-\lambda)(\lambda+\tfrac{\xi_{2}}{\Delta x})(\Delta x)^3.
\end{array}
\end{equation*}
Here we have used the relation $u^{(4)}(t,x) =u^{(4)}(0,x- c t)$.
\end{proof}
\begin{rem}
The consistency argument rests on the error estimates for the Hermite cubic polynomial interpolation, and the same argument can be applied to prove the consistency of the scheme \eqref{newu} and \eqref{tra}: These schemes are of fourth-order accuracy in space for $u$ and of the third-order accuracy in space for $u_x$. Note that $\Delta t>0$ is chosen arbitrary independent of the mesh size $\Delta x$. If one takes $\Delta t \sim \Delta x$, the schemes become third-order accurate scheme in time and space for $u$ and the second order in time and space for $u_x$.
\end{rem}

Let us denote the error in the numerical solution by $e^n_k=u^n_k -U^n_k $ and $r^n_k = v^n_k- V^n_k$.
From \eqref{u}, we see that the local errors $e^n_k$, $r^n_k$ satisfy the recursion relation
\begin{equation*}
\left(\begin{array}{c} e^{n+1}_k \\ \Delta x\; r^{n+1}_k
\end{array}\right)= A(\lambda)\left(\begin{array}{c} e^{n}_j \\
\Delta x\;{r}^{n}_j
\end{array}\right)+B(\lambda)\left(\begin{array}{c}  e^{n}_{j-1} \\ \Delta x \;{r}^{n}_{j-1} \end{array}\right)
+\left(\begin{array}{c} (E^{1})^n_j \\ \Delta x (E^{2})^n_j\end{array}
\right).
\end{equation*}

We multiply the equation by $e^{-i k\theta}$ and sum over $k\in \mathbb{Z}$; the outcome is
\begin{equation}\label{err}
\left(\begin{array}{c}\hat{e}^{n+1}(\theta) \\
\Delta x\;\hat{r}^{n+1}(\theta) \end{array} \right)=
 G_{\theta,\lambda}
\left(\begin{array}{c}\hat{e}^{n}(\theta) \\
\Delta x\;\hat{r}^{n}(\theta)\end{array}\right)+\hat{E}^n(\theta),
\end{equation}
where
\begin{equation*}% \label{onestep}
G_{\theta,\lambda}=A(\lambda)+\mbox{exp}(-i \theta)\,B(\lambda),\quad
\hat{E}^n(\theta) = \left(\begin{array}{c}
\sum_{k\in \mathbb{Z}} (E^{1})^n_{k-\ell} e^{-i k\theta}\\
\Delta x \sum_{k\in\mathbb{Z}} (E^{2})^n_{k-\ell} e^{-i k\theta}\end{array}
\right).
\end{equation*}
%and
%\[
%\hat{E}^n(\theta) = \left(\begin{array}{c}
%\sum_{k\in \mathbb{Z}} (E^{1})^n_{k-\ell} e^{-i k\theta}\\
%\Delta x \sum_{k\in\mathbb{Z}} (E^{2})^n_{k-\ell} e^{-i k\theta}\end{array}
%\right).
%\]

We investigate the stability of the map $G_{\theta,\lambda}$. Let us denote the eigenvalues of the amplification matrix $G_{\theta,\lambda}$ by $\rho_{1,\theta,\lambda}$ and $\rho_{2,\theta,\lambda}$  with $|\rho_{1,\theta,\lambda}|\le |\rho_{2,\theta,\lambda}|$.
\begin{thm}[Conditional stability]\label{thm:stability}
There exists constants $M_{\theta,\lambda}>0$ and $0 < \rho_{\theta,\lambda} \le 1$
depending on $(\theta,\lambda)\in [0,2\pi]\times[0,1]$, such that
\begin{equation}\label{boundG}
\|G_{\theta,\lambda}^n\|\le M_{\theta,\lambda} \rho^{n}_{\theta,\lambda}.
\end{equation}
%$\rho_{\theta,\lambda}=1$ if and only if $\theta=0,\; 2\pi$ or $\lambda=0,\;1$.
\end{thm}

The proof of the theorem rests on the following observations. We postpone their proofs in the Appendix.
\begin{lem}\label{lem:rhosimple}
For $0<\lambda < 1$ and $0\le \theta \le 2\pi$, the eigenvalues $\rho_{1,\theta,\lambda}$ and $\rho_{2,\theta,\lambda}$  are simple.
%Therefore, the amplification matrix $G_{\theta,\lambda}$ is diagonalizable.
\end{lem}
\begin{lem}\label{lem:rho_less_thanone}
$|\rho_{1,\theta,\lambda}| \le |\rho_{2,\theta,\lambda}| < 1$ for $(\theta,\lambda)\in (0,2\pi)\times(0,1)$.
\end{lem}

Now we give the proof of Theorem \ref{thm:stability}.
\begin{proof}
In case $\lambda=0$ or $\lambda=1$, $G_{\theta,\lambda}$ becomes trivial: identity map $I$ when $\lambda=0$, and $G_{\theta,1}=\exp(-i\theta)I$ when $\lambda=1$. Thus $\|G_{\theta,\lambda}^n\|=1$ for all $n\in\N$.
If $\theta=0, 2\pi$, the eigenvalues of $G_{\theta,\lambda}$
%$$
%G_{\theta,\lambda}=A(\lambda)+B(\lambda)\quad\mbox{and}\quad \lim_{n\to\infty}\,
%G_{\theta,\lambda}^n=\left(\begin{array}{cc} 1 & \frac{2\lambda-1}{6} \\ 0 &0 \end{array}\right).
%$$
are $\rho_{1,\theta,\lambda}=1-6\lambda + 6\lambda^2<1$ and $\rho_{2,\theta,\lambda}=1$, and we have
\[
G_{\theta,\lambda} = \left(
\begin{array}{cc}
1 & \frac{1-2\lambda}{6} \\
0 & 1 \\
\end{array}
\right)^{-1}
                     \left(
                       \begin{array}{cc}
                         1 & 0 \\
                         0 &  1-6\lambda + 6\lambda^2\\
                       \end{array}
                     \right)\left(
\begin{array}{cc}
1 & \frac{1-2\lambda}{6} \\
0 & 1 \\
\end{array}
\right).
\]
%where
%\[
%V=\left(
%\begin{array}{cc}
%1 & \frac{1-2\lambda}{6} \\
%0 & 1 \\
%\end{array}
%\right)
%\]
Thus we can take $\rho_{\theta,\lambda}=1$ and $M_{\theta,\lambda}= \|V^{-1}\|\|V\|$ where $V=\left(
\begin{array}{cc}
1 & \frac{1-2\lambda}{6} \\
0 & 1 \\
\end{array}
\right)$. %$\|G_{0,\lambda}^n\|$ is uniformly bounded in $n$.

Let us assume that $(\theta,\lambda)\in (0,2\pi)\times(0,1)$. The estimate \eqref{boundG} for $G_{\theta,\lambda}$ follows from Lemma \ref{lem:rhosimple} and Lemma \ref{lem:rhosimple}.
Indeed, by Lemma \ref{lem:rhosimple}, $G_{\theta,\lambda}$ is diagonalizable, i.e., there exists a invertible matrix $V_{\theta,\lambda}$
\begin{equation}\label{diagG}
  G_{\theta,\lambda}= V^{-1}_{\theta,\lambda}\left(
                         \begin{array}{cc}
                           \rho_{1,\theta,\lambda} & 0 \\
                           0 &\rho_{2,\theta,\lambda} \\
                         \end{array}
                       \right)V_{\theta,\lambda},
\end{equation}
and thus
%\begin{equation*}
$\|G^n_{\theta,\lambda}\|\le \|V^{-1}_{\theta,\lambda}\|\|V_{\theta,\lambda}\||\rho_{2,\theta,\lambda}|^n$
%\end{equation*}
for all $n \in \N$. Consequently, we have \eqref{boundG} by setting $M_{\theta,\rho}:=\|V^{-1}_{\theta,\lambda}\|\|V_{\theta,\lambda}\|$ and $\rho_{\theta,\lambda}:=|\rho_{2,\theta,\lambda}|$. The inequality $\rho_{\theta,\lambda} < 1$ follows from Lemma \eqref{lem:rhosimple}.
\end{proof}
\begin{thm}
\begin{equation}\label{errorest}
\left\|\left(\begin{array}{c}\hat{e}^{n}(\theta) \\
\Delta x\; \hat{r}^{n}(\theta) \end{array} \right)\right\| \le
%\le M_\theta(\lambda) rho_k^n
%\left\|\left(\begin{array}{c}\hat{e}^0_k \\
%\hat{r}^{0}_k \end{array}
%\right)\right\|+
%(\Delta x)^4M_{\theta,\lambda}\left(\sum_{i=0}^{n-1}\rho^i_{\theta,\lambda}\right)|u^{(4)}(0,\cdot)|_{L^\infty(\mathbb{R})}.
M_{\theta,\lambda} \rho^{n}_{\theta,\lambda}\left\|\left(\begin{array}{c}\hat{e}^{0}(\theta) \\
\Delta x\; \hat{r}^{0}(\theta) \end{array} \right)\right\| +
M_{\theta,\lambda}\left(\sum_{i=0}^{n-1}\rho^i_{\theta,\lambda}\right)\|\hat{E}^n(\theta)\|
%\sim M_{\theta,\lambda}(\Delta x)^4\left(\sum_{i=0}^{n-1}\rho^i_{\theta,\lambda}\right).
%(1-\rho^n_{\theta,\lambda})}{1-\rho^{}_{\theta,\lambda}}|u^{4}|_{L^\infty(\mathbb{R})}.
\end{equation}
for all $n\in \N$.
\end{thm}
\begin{proof}
Apply \eqref{err} repeatedly and remember \eqref{diagG}:
\begin{align*}
\left(
     \begin{array}{c}
      \hat{e}^{n}(\theta) \\
      \Delta x\; \hat{{r}}^{n}(\theta) \\
     \end{array}
   \right)
  & =
   G^{n}_{\theta,\lambda}\left(
     \begin{array}{c}
      \hat{e}^{0}(\theta) \\
     \Delta x\;  \hat{{r}}^{0}(\theta) \\
     \end{array}
   \right)
+(G^{n-1}_{\theta,\lambda} + \cdots + G_{\theta,\lambda} + I)
\hat{E}^n(\theta)\\
   &
  = G^{n}_{\theta,\lambda}\left(
     \begin{array}{c}
      \hat{e}^{0}(\theta) \\
     \Delta x\;  \hat{{r}}^{0}(\theta) \\
     \end{array}
   \right) +  V^{-1}_{\theta,\lambda}\left(
                         \begin{array}{cc}
                           \sum_{i=0}^{n-1}\rho^{i}_{1,\theta,\lambda} & 0 \\
                           0 &\sum_{i=0}^{n-1}\rho^{i}_{2,\theta,\lambda} \\
                         \end{array}
                       \right)V_{\theta,\lambda}\hat{E}^n(\theta),
\end{align*}
and so \eqref{errorest} follows from \eqref{boundG}.
\end{proof}

It is easy to see that the residual error is of the order $(\Delta x)^4$, i.e.,
\begin{equation*}
\|\hat{E}^n(\theta)\| \sim
(\Delta x)^4|u^{(4)}(0,\cdot)|_{L^\infty(\R)}.
\end{equation*}
Let us assume that the initial condition for the derivative satisfies $|v^{0}_k - u_x(0,x_k)| \sim (\Delta x)^3$. Then \eqref{errorest} derives an error estimate: Let us denote the finial time by $T$ and the total number of iteration by $n$.
\begin{cor}
\begin{equation*}
|\hat{e}^{n}(\theta)| \sim \frac{T}{\Delta t}(\Delta x)^4,\quad
|\hat{r}^{n}(\theta)|\sim \frac{T}{\Delta t}(\Delta x)^3
\end{equation*}
for $\theta \in [0,2\pi]$.
\end{cor}
\begin{rem} We have observed through numerical computation that
\begin{equation*}
M_{\theta,\lambda} = \|V^{-1}_{\theta,\lambda}\|\|V_{\theta,\lambda}\| \le  3.6453,
\end{equation*}
uniformly for $(\theta,\lambda)\in[0,2\pi]\times[0,1]$.
Thus, we conjecture the constant $M_{\theta,\lambda}$ appearing in \eqref{boundG} and \eqref{errorest} is uniformly
bounded in $\theta$ and $\lambda$. The uniform boundness of $M_{\theta,\lambda}$ leads to a rigorous proof of the error estimates
\begin{align*}
\|u(t_n,\cdot)-u^n\|_{\Delta x}  \sim \frac{T}{\Delta t}(\Delta x)^4, \quad
\|u_x(t_n,\cdot)-v^n\|_{\Delta x}  \sim \frac{T}{\Delta t}(\Delta x)^3.
\end{align*}
We reserve the proof of the boundedness of $M_{\theta,\lambda}$ for a future work.
\end{rem}

\subsection{Numerical results}
We report numerical results for the advection equation with variable coefficient. Let us consider the advection equation $u_t + (c(x)u)_x=0$, $c(x)=\frac{1}{\cos(4\pi x) + 2}$ for $x\in[0,1], \; t>0$, with periodic boundary condition.
As an initial condition, we take $u(0,x)=\exp(-\tfrac{(x-0.2)^2}{0.05^2})$, which can be regarded as periodic in practice. One can easily see that the solution at time $t=2, 4, 6, \ldots$ is given by the initial condition $u(0,x)$, i.e., $u(2m,x)=u(0,x)$ for $m\in\N$.
%We compare a numerical solution at time $t=2m$ to the initial condition, which is the exact solution at the time. \\
%and the advection equation $\partial_t u + c(x)u_x=0$

We report the accuracy of the CIP method. In this numerical test, the time step size is fixed to be $\Delta t= 0.1$ for each mesh size $N^{-1}$, $N\in\{50, 100, 200, 400, 800, 1600\}$.
The numerical solutions at time $t=2$ are computed and compared to the exact solution. The number of time integration is $2/\Delta t = 20$ for all $N$. For each mesh size, the error in the numerical solutions is measured by $\ell^1$, $\ell^2$ and $\ell^\infty$ norm:
\begin{equation}\label{numericalerror}
\epsilon_\infty =\max_k \max_{x\in[0,1]}|u^n_k - u(t,x_k)|,\quad  \epsilon_i = \frac{|u^n - u(t,\cdot)|_{\ell^i}}{|u(t,\cdot)|_{\ell^i}}, \quad i=1,\;2.
\end{equation}
where $u^n=\{u^n_k\}_{k=1}^N$ is the numerical solution and $u(t,\cdot)=\{u(t,x_k)\}_{k=1}^N$ is the exact solution at the grids. The characteristic equation,
$ \frac{d x(s)}{d s} = c(x(s)), x(\Delta t ) = x_k$, is solved backward in time to find the location $y_k=x(0)$ for each $k$. We employ Matlab build in function \textsf{ode23} to solve the characteristic equation.
%with \textsf{RelTol} = \textsf{1e-8}, \textsf{AbsTol} = \textsf{1e-10} for the numerical integration of
%Plot of Figure \ref{fig:snapshotsCIPLargeadv} shows four snapshots of the numerical solutions ( circle ) and the velocity field $c(x)$ (solid line).

\begin{figure}[h]
\centering
\includegraphics[width=6cm, bb=100 230 500 540]{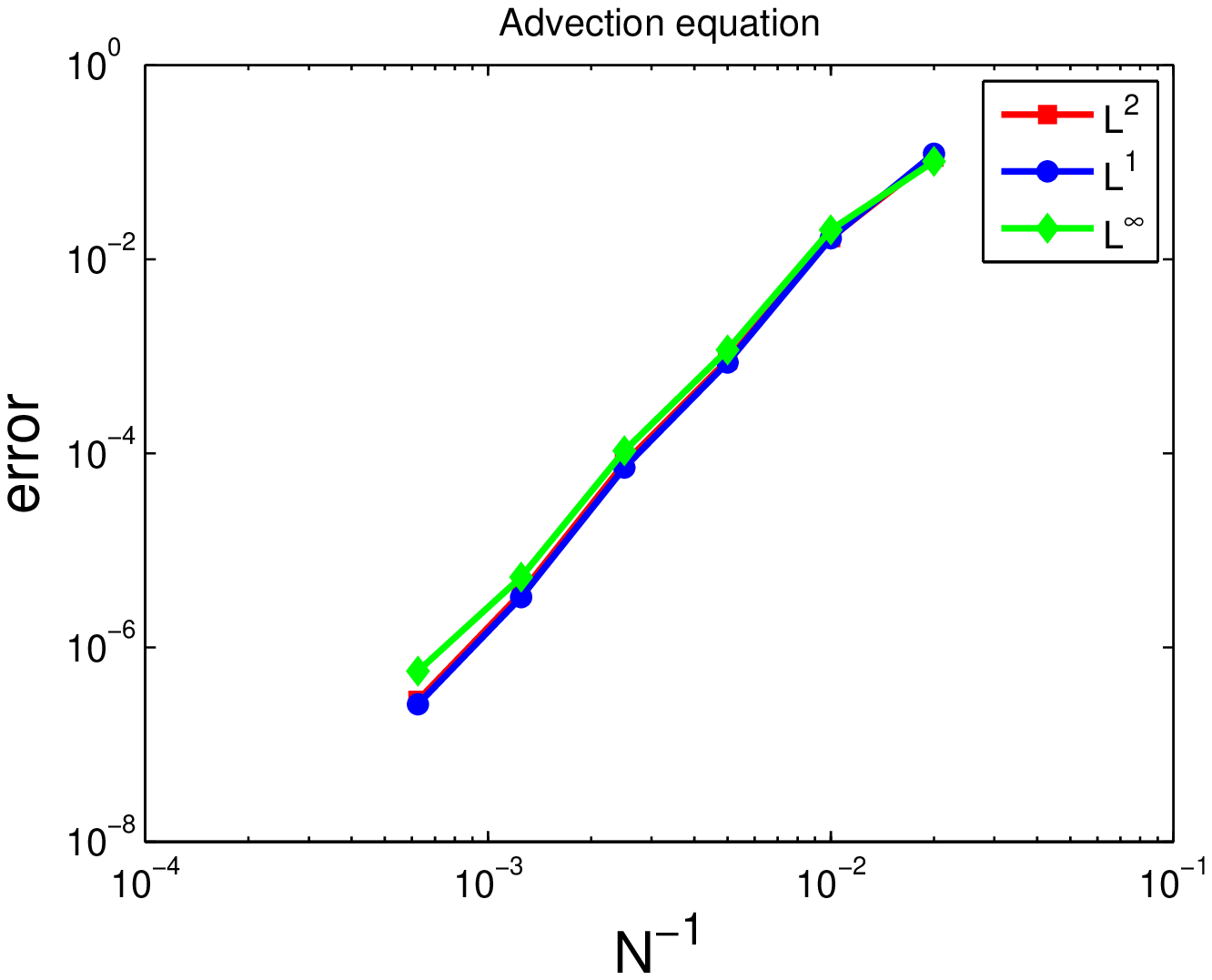}
\includegraphics[width=6cm, bb=100  230 500 540]{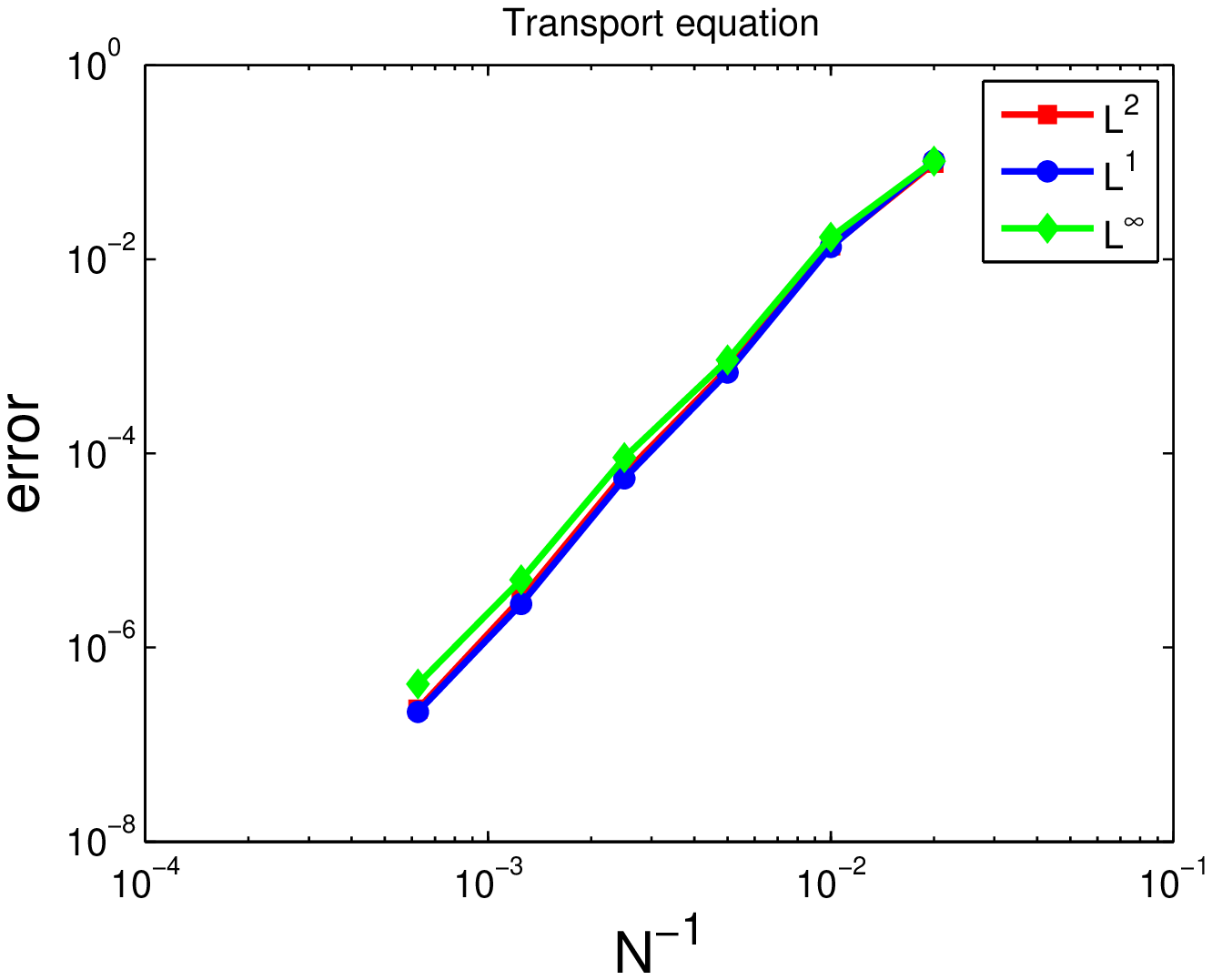}
\caption{The fourth order convergence in space against mesh size $N^{-1}$ of the numerical error in the numerical solution by the CIP. The error is computed by \eqref{numericalerror} at $t=2$.}
\label{fig:errorCIPLarge}
\end{figure}

Plot of Figure \ref{fig:errorCIPLarge} (left) and Table \ref{table:traadv} (left) show the fourth order convergence of the method. We also report the performance of the CIP for the transport equation $u_t + c(x)u_x=0$, $c(x)=\frac{1}{\cos(4\pi x) + 2}$ for $x\in[0,1], \; t>0$, with periodic boundary condition in Figure \ref{fig:errorCIPLarge} (right) and Table \ref{table:traadv} (right).
%
%\begin{figure}[h]
%\centering
%\includegraphics[width=0.5\textwidth, bb=100 230 500 540]{CIPLargeCFLadv.eps}
%\caption{The snapshots of the numerical solutions at time $t\in\{0,\; 0.2,\; 0.8,\; 1.4\}$ of the advection equation ( blue circle ). The velocity $c(x) = \frac{1}{\cos4\pi x + 2}$ is depicted in red solid line.}
%\label{fig:snapshotsCIPLargeadv}
%\end{figure}

%\begin{table}[h]
%\centering\caption{Numerical results for the advection equation.}
%\begin{tabular}{|c|c|c|c|c|c|c|}
%  \hline
%  % after \\: \hline or \cline{col1-col2} \cline{col3-col4} ...
%  $N$ & 50 & 100 & 200 & 400 & 800 & 1600 \\
%   \hline
%  $\epsilon_1$  &1.03e-1 & 1.33e-2 & 6.79e-4&  5.50e-5&  2.80e-6 & 2.16e-7  \\
%  $\epsilon_2$& 9.75e-2 & 1.37e-2 & 7.43e-4&  6.17e-5 & 3.19e-6&  2.31e-7   \\
%  $\epsilon_\infty$ & 1.02e-1 & 1.68e-2 & 9.18e-4 & 9.04e-5 & 4.97e-6  &4.19e-7  \\
%  \hline
%\end{tabular}\label{table:adv}
%\end{table}
%\begin{table}[h]
%\centering\caption{Numerical results for the transport equation.}
%\begin{tabular}{|c|c|c|c|c|c|c|}
%  \hline
%  % after \\: \hline or \cline{col1-col2} \cline{col3-col4} ...
%  $N$ & 50 & 100 & 200 & 400 & 800 & 1600 \\
%   \hline
%  $\epsilon_1$  & 1.22e-1 & 1.64e-2 & 8.62e-4  &7.11e-5  &3.30e-6  &2.60e-7  \\
%  $\epsilon_2$  &1.12e-1  &1.66e-2  &9.13e-4 & 7.92e-5 & 3.63e-6 & 2.85e-7  \\
%  $\epsilon_\infty$  &1.01e-1 & 2.01e-2 & 1.16e-3 & 1.06e-4 & 5.30e-6 & 5.72e-7  \\
%  \hline
%\end{tabular}\label{table:tra}
%\end{table}

\begin{table}[h]\scriptsize
\centering\caption{Numerical errors in the solutions for the advection and the transport equation.}
\begin{tabular}{c|cccccc|cccccc|}
  &\multicolumn{6}{|c|}{Advection equation}& \multicolumn{6}{c|}{Transport equation} \\
  \hline
  $N$ & \scriptsize{50} & 100 & 200 & 400 & 800 & 1600  & 50 & 100 & 200 & 400 & 800 & 1600 \\
   \hline
  $\epsilon_1$  &1.03e-1 & 1.33e-2 & 6.79e-4&  5.50e-5&  2.80e-6 & 2.16e-7 & 1.22e-1 & 1.64e-2 & 8.62e-4  &7.11e-5  &3.30e-6  &2.60e-7 \\
  $\epsilon_2$& 9.75e-2 & 1.37e-2 & 7.43e-4&  6.17e-5 & 3.19e-6&  2.31e-7  &1.12e-1  &1.66e-2  &9.13e-4 & 7.92e-5 & 3.63e-6 & 2.85e-7 \\
  $\epsilon_\infty$ & 1.02e-1 & 1.68e-2 & 9.18e-4 & 9.04e-5 & 4.97e-6  &4.19e-7  &1.01e-1 & 2.01e-2 & 1.16e-3 & 1.06e-4 & 5.30e-6 & 5.72e-7 \\
  \hline
  \end{tabular}\label{table:traadv}
\end{table}

\section{Immersed interface method for CIP}\label{sec:IIMCIP}
%(Matlab code:\verb"CIP1dCubic.m", \verb"IIMCIP1dCubic.m", \verb"IIMCIP1d.m")
In this section we develop the immersed interface method for CIP
(IIM-CIP for short) for transport equations with discontinues
coefficient;
\begin{equation}\label{equ:discon_tra}
u_t + c(x)u_x = 0, \quad t>0, \ x\in \R,\qquad u(0,x)=u_0(x),\;\;
x\in \R,
\end{equation}
where $c>0$ is a piecewise smooth that has a jump across $x=\alpha$.
We assume that the interface is located in an
interval $[x_{j-1},x_j]$. That is, $c=c^-$ in $x<\alpha$ and $c=c^+$ in
$x>\alpha$ with $c^-\neq c^+$. We call the interval $[x_{j-1},x_j]$
\textit{an irregular interval} or \textit{an irregular cell}, and
the point $x_j$ \textit{an irregular point}. The interface condition
on $u$ is imposed according to the physical phenomena under
consideration. We take $[u]=0$ as an interface condition at the interface $\alpha$, where $[u]:= \lim_{x\to \alpha^+}u^+(x) - \lim_{x\to
\alpha^-}u^-(x)$. Here $u^-(x) = u(x)|_{[x_{j-1},\alpha]}$ and
$u^+(x) = u(x)|_{[\alpha,x_j]}$. The treatment of an interface
condition $[cu]=0$ will be briefly discussed at the end of this
section.

From the interface condition coupled with the equation,
we know the flux $(cu)_x=cu_x$ is also continuous across the interface,
and thus the derivative $u_x$ has a jump discontinuity at the interface. Because of the discontinuity in $u_x$,
the standard CIP using a single profile in an interval
will not provide an accurate solution.

We begin with the construction of a piecewise Hermite cubic polynomial that approximates to the solution $u(t_n,x)$ in the interval.
Next, we derive a time integration formula for the exact solution and its derivative.
Lastly, we propose a numerical scheme (IIM-CIP scheme) to update the numerical solution $u^{n}_j$ and the derivative $v^{n}_j$ to the next time level.

\textsf{Piecewise cubic polynomial.} Let us consider to approximate the solution $u(t_n,x)$, $x\in [x_{j-1},x_j]$. The interface condition $[u]=0$ coupled with the equation \eqref{equ:discon_tra} yields the relations $[c^{k}\frac{\partial^ku}{\partial x^k}]=0$, $k\in\N$. See \cite{Zhang+LeVeque-immeintemethacou:97} for the derivation.  Obviously, it is impossible to construct a single polynomial interpolation that satisfies the relations.
We introduce two polynomials of the form $H^{\pm}(x) = \sum_{k=0}^3 \frac{a^{\pm}_k}{k!}(x-\alpha)^k$.
We determine the eight unknowns via the interface relations and the interpolation conditions.
%Let us denote a piecewise cubic polynomial $H(x)$ by $H|_{[x_{j-1},\alpha]}=H^-|_{[x_{j-1},\alpha]}$, $H|_{[\alpha,x_{j}]}=H^+|_{[\alpha,x_{j}]}$. Then we impose the conditions;
\begin{align}
&\left\{\begin{array}{cc}
  H^-(\alpha)=H^+(\alpha),                         & c^-H^-_x(\alpha)= c^+H^+_x(\alpha), \\
  (c^-)^2H^-_{xx}(\alpha)=(c^+)^2H^+_{xx}(\alpha), & (c^-)^3H^-_{xxx}(\alpha)=(c^+)^3H^+_{xxx}(\alpha),
\end{array}\right. \label{IIM}    \\
 \nonumber \\
&  H^-(x_{j-1})=u^n_{j-1},\quad H_x^-(x_{j-1})=v^n_{j-1},\quad H^+(x_{j})=u^n_{j},\quad H_x^+(x_{j})=v^n_{j}.\label{interp1}
\end{align}
The first equations yield
\begin{align}
a^{-}_0 =
a^{+}_0,\quad
c^-a^{-}_1=
c^+a_1^+, \quad
(c^-)^2 a^{-}_2=
(c^+)^2 a^{+}_2, \quad
(c^-)^3 a_3^{-} =
(c^+)^3 a_3^{+},\label{interface1}
\end{align}
and thus, introducing new parameters
$a=(a_0,\ldots,a_3)^\top$, we can write
$
H^\pm(x)=\sum_{\ell=0}^3 \frac{a_\ell}{\ell !}\left(\frac{x-\alpha}{c^\pm\Delta x}\right)^\ell.
$ Then using the interpolation conditions \eqref{interp1}, we obtain the system $A a = f$, where
\begin{equation*}%\label{A}
  A=\left(
      \begin{array}{cccc}
        1 & \frac{\theta}{c^+} & \frac{\theta^2}{2(c^+)^2} & \frac{\theta^3}{3!(c^+)^3} \\
        0 & 1 & \frac{\theta}{c^+} & \frac{\theta^2}{2(c^+)^2} \\
        1 & \frac{(\theta-1)}{c^-} & \frac{(\theta-1)^2}{2(c^-)^2} & \frac{(\theta-1)^3}{3!(c^-)^3} \\
        0 & 1 & \frac{(\theta-1)}{c^-} & \frac{(\theta-1)^2}{2(c^-)^2}
      \end{array}
    \right),\quad
    f=\left(
        \begin{array}{c}
          u^n_{j} \\
         \Delta x\; v^n_{j}\\
          u^n_{j-1} \\
         \Delta x\; v^n_{j-1}
        \end{array}
      \right),
\end{equation*}
and $\theta = \frac{x_j-\alpha}{\Delta x}$. The determinant of $A$,
$\det A = \frac{(c^-\theta+ c^+(1-\theta))^4}{12(c^+c^-)^4}$, is positive for all $0\le \theta \le 1$ and thus the coefficient $a$ is uniquely determined.
We define a piecewise polynomial $H(x)$ on $[x_{j-1},x_j]$ by $H|_{[x_{j-1},\alpha]}=H^-|_{[x_{j-1},\alpha]}$, $H|_{[\alpha,x_{j}]}=H^+|_{[\alpha,x_{j}]}$ and call it \textit{the immersed interface cubic polynomial} to the data set $f$. We also denote the piecewise polynomial by $H^{\pm}$.
%\begin{equation}\label{A}
%  A=\left[
%      \begin{array}{cccc}
%        1 & \frac{d_+}{c_+} & \frac{d_+^2}{2c_+^2} & \frac{d_+^3}{3!c_+^3} \\ \\
%        0 & \frac{1}{c_+} & \frac{d_+}{c_+^2} & \frac{d_+^2}{2c_+^3} \\ \\
%        1 & \frac{d_-}{c_-} & \frac{d_-^2}{2c_-^2} & \frac{d_-^3}{3!c_-^3} \\ \\
%        0 & \frac{1}{c_-} & \frac{d_-}{c_-^2} & \frac{d_-^2}{2c_-^3} \\
%      \end{array}
%    \right],\quad
%    f=\left[
%        \begin{array}{c}
%          u^n_{j} \\
%          v^n_{j} \\
%          u^n_{j-1} \\
%          v^n_{j-1} \\
%        \end{array}
%      \right].
%\end{equation}
%where $d_{-}=x_{j-1}-\alpha$ and $d_+=x_j-\alpha$.
\begin{figure}
  \centering
  \subfloat[{}]{\label{IIMprofile1d}\includegraphics[width=0.5\textwidth]{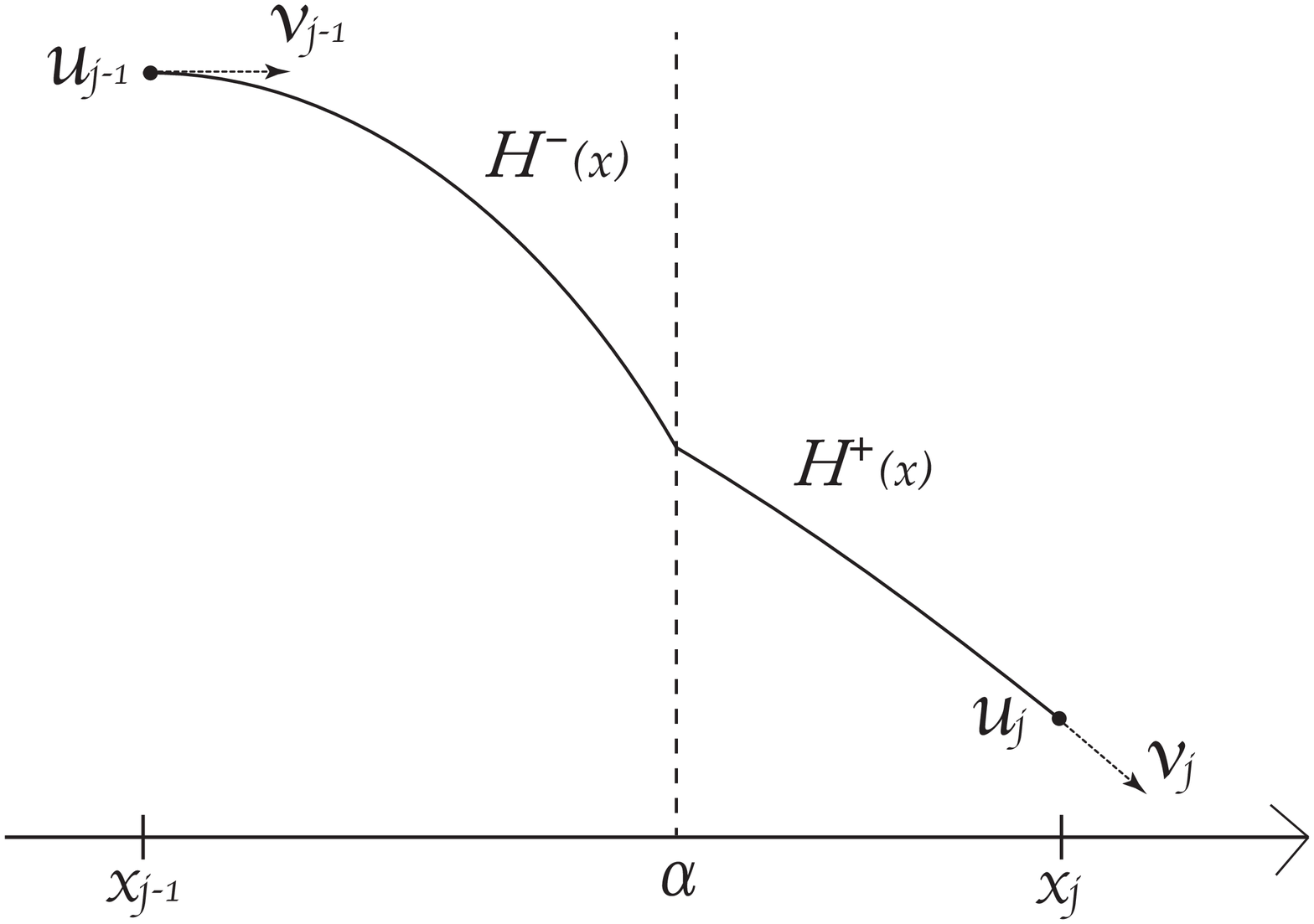}}
  \subfloat[{}]{\label{IIMprofileExtUpdate}\includegraphics[width=0.5\textwidth]{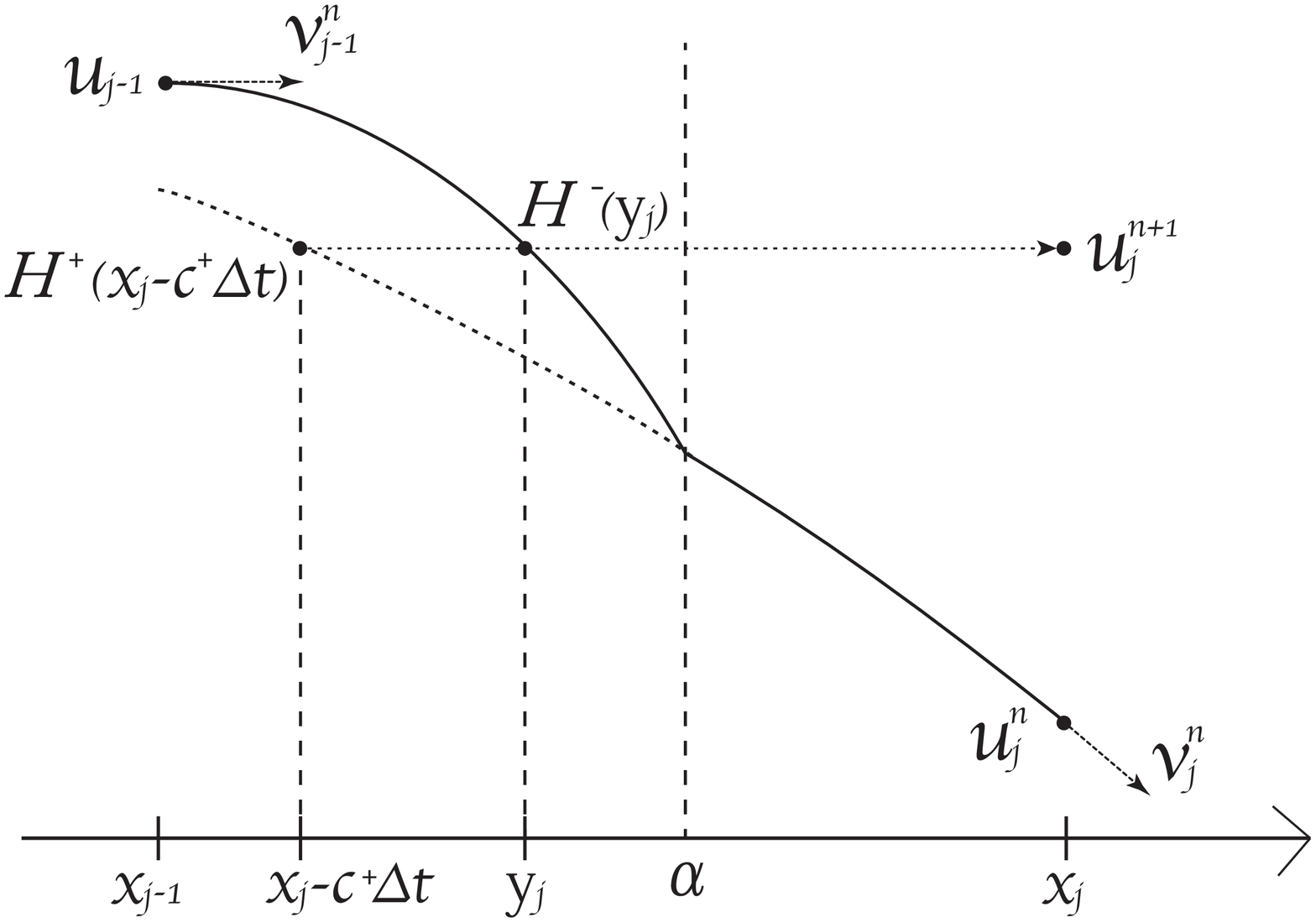}}
  \caption{(a) The piecewise cubic polynomial $H^\pm$ on the irregular interval. It is continuous at the interface $\alpha$ and has jump discontinuity in the one-sided derivative at the interface.
  (b) Graphical illustration of the IIM-CIP scheme. The value $H^+(x_j-c^+\Delta t)$ can be used for updating the numerical solution at the time level $t_{n+1}$ at the grid $x_j$.}
\end{figure}

Figure \ref{IIMprofile1d} illustrates the immersed interface cubic polynomial $H^\pm(x)$ on the interval $[x_{j-1},x_j]$.
One can see that the immersed interface cubic polynomial is continuous at the interface but has discontinuity in the one-sided derivative at the point.

We have observed in Section \ref{sec:error} that the accuracy of the Hermite cubic polynomial is of the fourth order in the function value
and is of third order in its derivative and directly affects to the accuracy in the one step map of the CIP method. As for the accuracy of the immersed interface cubic polynomial, we have the following:
%The piecewise cubic interpolation $H(x)$ is indeed an accurate approximation to $u(x)$.
Let us consider the immersed interface cubic polynomial $h^\pm(x)$ to the exact solution $u(t_n,x)$.
From the interface relations \eqref{IIM}, the polynomials $h^+$ and $h^-$ can be written in the form
$ h^+(x)=\sum_{\ell=0}^3 \frac{a_\ell}{\ell !}\left(\frac{x-\alpha}{c^+\Delta x}\right)^\ell$ and $ h^-(x)=\sum_{\ell=0}^3 \frac{a_\ell}{\ell !}\left(\frac{x-\alpha}{c^-\Delta x}\right)^\ell$.
The interpolation condition with exact solution as a data set,
\begin{equation*}
%&[h]=0,\quad [ch_x]=0,\quad [c^2h_{xx}]=0,\quad [c^3h_{xxx} ]=0,   \\
h^-(x_{j-1})=u(t_n,x_{j-1}),\   h_x^-(x_{j-1})=u_x(t_n,x_{j-1}),\ h^+(x_{j})=u(t_n,x_{j}),\ h_x^+(x_{j})=u_x(t_n,x_{j}),
\end{equation*}
leads to the system $A a = f$ for the coefficient $a$,
%The coefficient $a$ is given by $A a = f$,
where
\[
f=(u^+(t_n,x_j),\;\Delta x\; u_x^+(t_n,x),\; u^-(t_n,x_{j-1}),\; \Delta x\; u_x^-(t_n,x))^\top.
\]
\begin{thm}[Accuracy of the immersed interface cubic polynomial]\label{thm:IIM}
Let $u(t,x)$ be the solution of \eqref{tra}. We have the following error estimates
\begin{equation*}
  u(t_n,x) - h(x) = O((\Delta x)^4),\quad  u_x(t_n,x) - h_x(x) = O((\Delta x)^3),  \quad x\in[x_{j-1},x_j].
\end{equation*}
Here the one-sided derivative is taken at the interface $\alpha$.
\end{thm}
\begin{proof}
The Taylor expansion for $u^\pm(t_n,x)$ at the interface $\alpha$ and the interface relations
\begin{equation*}
  \begin{array}{l}
u^-(t_n,\alpha)=u^+(t_n,\alpha), \quad c^{-}u^-_{x}(t_n,\alpha)=c^{+}u^+_{x}(t_n,\alpha),\\
(c^{-})^2u^-_{xx}(t_n,\alpha)=(c^{+})^2u^+_{xx}(t_n,\alpha), \quad
(c^{-})^3u^-_{xxx}(t_n,\alpha)=(c^{+})^3u^+_{xxx}(t_n,\alpha),
  \end{array}
\end{equation*}
lead to
%\begin{equation}
%    f=Ab+
%\end{equation}
% where $b$  $b_0=u^+(t_n,\alpha)$, $b_1=c_{+}u^+_{x}(t_n,\alpha)$,
%$b_2=c_{+}^2u^+_{xx}(t_n,\alpha)$ and $b_3=c_{+}^3u^+_{xxx}(t_n,\alpha)$. Thus
%$A(a-b)=r$.
\begin{equation*}
\begin{array}{lll}
u^+(t_n,x_j)
&=& b_0 + b_1\frac{\theta}{c^+}+\frac{b_2}{2}\left(\frac{\theta}{c^+}\right)^2
+ \frac{b_3}{3!}\left(\frac{\theta}{c^+}\right)^3 + O((\Delta x)^4),\\
\Delta x\; u_x^+(t_n,x)
&=&\frac{b_1}{c^+}+b_2\frac{\theta}{(c^+)^2}+\frac{b_3}{2}\frac{\theta^2}{ (c^+)^3} + O((\Delta x)^4),\\
u^-(t_n,x_{j-1})
&=& b_0 + b_1\frac{\theta-1}{c^-}+\frac{b_2}{2}\left(\frac{\theta-1}{c^-}\right)^2
+ \frac{b_3}{3!}\left(\frac{\theta-1}{c^-}\right)^3 + O((\Delta x)^4),\\
\Delta x\; u_x^-(t_n,x)
&=&\frac{b_1}{c^-}+b_2\frac{\theta-1}{( c^-)^2}+\frac{b_3}{2}\frac{(\theta-1)^2}{( c^-)^3} + O((\Delta x)^4),
\end{array}
\end{equation*}
where $ b_0=u^+(t_n,\alpha)$, $b_1=c^{+}u^+_{x}(t_n,\alpha)$,
$b_2=(c^{+})^2u^+_{xx}(t_n,\alpha)$ and $b_3=(c^{+})^3u^+_{xxx}(t_n,\alpha)$. Thus we have
$A (a-b) =(r_1,r_2,r_3,r_4)^\top$ where $r_i = O((\Delta x)^4)$ for $i=1,2,3,4$.

Since the components of $A$ is of order $O(1)$, the components of the inverse $A^{-1}$ is also of order $O(1)$.
Therefore, the equation $A (a- b) = r$ implies that $
  a_k-b_k = O((\Delta x)^4)$ for $k=0,1,2,3.$ Thus we obtain the desired estimates.
\end{proof}

\noindent
\textsf{Update formula.}
Now that we have constructed the immersed interface cubic polynomial in the irregular interval, we attempt to build an
update formula for numerical solutions with the
aid of the method of characteristic. Throughout of this section we assume that CFL number is less or equal to 1, i.e., $c^\pm \Delta t \le \Delta x.$ Thus $y_k$ is always included in the interval $[x_{k-1},x_k]$ for all $k$.

If $x_j$ is an regular point, the numerical solution is updated by \eqref{u}
 with $\lambda=\frac{c^- \Delta t}{\Delta x}$ or $\lambda=\frac{c^+ \Delta t}{\Delta x}$
 depending on the velocity at the grid $x=x_j$.
 Let us assume that $x_j$ is an irregular grid,
 and let us consider the characteristic curve $x=x(s)$ with $x(\Delta t)=x_j$.
 Let us denote the upwind point $x(0)$ by $y_j$. %, which we
%will determine below, in the interval $[x_{j-1},x_j]$ at $s=0$,
%arriving at $x_j$ at $s=\Delta t$.
\begin{Prop}
Suppose that there exists an interface in the interval $[x_{j-1},x_j]$.
%The upwind point $y_k$, the solution of $u^+(t_{n+1},x_k)$ and its derivative $u^+_x(t_{n+1},x_k)$ at the irregular grid are given by,
If $\alpha \le x_j -c^+\Delta t$,
\begin{equation*}\begin{array}{l}
y_j = x_j -c^+\Delta t, \\
u^+(t_{n+1},x_j)=u^+(t_n,y_j ),\quad u^+_x(t_{n+1},x_j)=u^+_x(t_n,y_j).
\end{array}
\end{equation*}
If $x_j -c^+\Delta t\le \alpha$,
\begin{equation*}
\begin{array}{l}
y_j = \alpha + \frac{c^-}{c^+}(x_j-\alpha)-c^-\Delta t,\\
u^+(t_{n+1},x_j)= u^-(t_n,y_j),\quad u_x^+(t_{n+1},x_j)=\frac{c^-}{c^+}u_x^-(t_n,y_j).
\end{array}
\end{equation*}
\end{Prop}
\begin{proof}
Let us assume that $\alpha \le x_j -c^+\Delta t$.
%There arise two cases: (a) $\alpha \le y_k=x(0)$, (b) $ y_k \le \alpha$.
% We first deal with the case (a).
 The characteristic ODEs
\eqref{ode:x}--\eqref{ode:p1} become $\dot{x}(s)=c^+$,
$\dot{z}(s)=\dot{p_1}(s)=\dot{p_2}(s)=0$, and thus we have $y_j=x_j-c^+\Delta
t$,  $u^+(t_{n+1},x_j)=u^+(t_n,x_j-c^+\Delta t)$ and
$u^+_x(t_{n+1},x_j)=u^+_x(t_n,x_j-c^+\Delta t)$.

Next let us assume that $x_j -c^+\Delta t \le \alpha$.
There exists $s^*>0$ such that $\alpha=x(s^*)$.
The characteristic curve obeys
\begin{equation*}
  \dot{x}(s) =
      c^-, \quad ( 0\le s \le s^*), \quad
      \dot{x}(s) =
      c^+, \quad  (s^*\le s\le \Delta t).
\end{equation*}
It is obvious to see that
$
     \alpha-y_j = c^- s^*, \quad
      x_j - \alpha= c^+(\Delta t - s^*),
$
which yields
$
  y_j = \alpha + \frac{c^-}{c^+}(x_j-\alpha)-c^-\Delta t, \quad
  s^* = \Delta t- \frac{x_j -\alpha}{c^+}.
$
The characteristic ODE for $z(s)=u(s,x(s))$ becomes
$\frac{d}{ds}{u^-(s,x(s))}=0, ( 0\le s \le s^*)$ and
$\frac{d}{ds}{u^+(s,x(s))}=0, ( s^*\le s \le \Delta t)$. Using the
interface condition $[u]=0$, it follows
$u^+(t_{n+1},x_j)=u^+(s^*,\alpha) = u^-(s^*,\alpha) = u^-(t_n,y_j)$.
Similarly, with the characteristic ODE for $u_x(s,x(s))$ and the
relation $[c u_x]=0$, we obtain $u_x^+(t_{n+1},x_j)=u_x^+(s^*,\alpha) =
\frac{c^-}{c^+}u_x^-(s^*,\alpha) =  \frac{c^-}{c^+}u_x^-(t_n,y_j)$.
\end{proof}

Thus we arrive at a CIP scheme for the equation \eqref{equ:discon_tra}: Let $H^+$ and $H^-$ be cubic polynomials defined by \eqref{IIM} and \eqref{interp1}.
If $\alpha\le x_j -c^+\Delta t$,
\begin{equation}\label{new1}
\begin{array}{l}
u^{n+1}_j =H^+(x_j-c^+\Delta t),\quad v^{n+1}_j= H^+_x(x_j-c^+\Delta t).
\end{array}
\end{equation}
and, if $x_j -c^+\Delta t\le \alpha$,
\begin{equation}\label{new2}
\begin{array}{l}
u^{n+1}_j =H^-(y_j),\quad v^{n+1}_j= \dfrac{c^-}{c^+}H_x^-(y_j).
\end{array}
\end{equation}
where $y_j = \alpha + \frac{c^-}{c^+}(x_j-\alpha)-c^-\Delta t$. \\
It seems one must appropriately choose either \eqref{new1} or \eqref{new2} according to the condition $\alpha\le x_j -c^+\Delta t$ or $x_j -c^+\Delta t\le \alpha$. However, no need for this special treatment arises and the update scheme at the irregular point is solely based on \eqref{new1}. Indeed we can show that if $x_j -c^+\Delta t\le  \alpha $, the numerical solution is also given by \eqref{new1}, i.e.,
\begin{equation}\label{Hm=Hp}
u^{n+1}=H^-(y_j) =  H^+(x_j - c^+\Delta t), \quad  v^{n+1}=\frac{c^-}{c^+}H_x^-(y_j)
= H^+_x(x_j - c^+\Delta t).
\end{equation}
From $\frac{c^{+}}{c^-}(y_j-\alpha) = x_j -c^+\Delta t -\alpha$ and \eqref{interface1}, we have
\begin{align*}
&a^{-}_\ell(y_j-\alpha)^\ell= a^{+}_\ell\left(\frac{c^{+}}{c^-}\right)^\ell(y_j-\alpha)^\ell
= a^{+}_\ell(x_j-c^+\Delta t-\alpha)^\ell,\\
&\frac{c^-}{c^+}a^{-}_\ell(y_j-\alpha)^{\ell-1}
=a^{+}_\ell\left(\frac{c^{+}}{c^-}\right)^{\ell-1}(y_j-\alpha)^{\ell-1} =
a^{+}_\ell(x_j-c^+\Delta t-\alpha)^{\ell-1},
\end{align*}
which implies \eqref{Hm=Hp}. \\
%Therefore the numerical solution is also given by \eqref{new1}.
Figure \ref{IIMprofileExtUpdate} illustrates that we can use $H^+(x_j-c^+\Delta t)$ and its derivative to update the numerical solutions $u^{n+1}_j$ and $v^{n+1}_j$ at the irregular grid even when $x_j-c^+\Delta t\le \alpha$ and \eqref{new2} is not necessary for the computation of the numerical solution. \\
Here is the summary of the IIM-CIP.
\begin{description}
%\item[IIM-CIP0] For each $k$, find the interval $[x_{j-1},x_j]$ that includes $y_j$.
\item[IIM-CIP1] Construct the cubic profile on each interval. If the interval is irregular, then construct the immersed interface cubic polynomial.
\item[IIM-CIP2] Update $u^{n+1}_j$ and $v^{n+1}_j$: If $x_j$ is a regular grid, use \eqref{u} with $\lambda=\frac{c\Delta t}{\Delta x}$ where $c=c^-$ if $x_j<\alpha$ or $c=c^+$ if $\alpha< x_j$. If $x_j$ is the irregular grid, use \eqref{new1}.
\end{description}

It is worth pointing out that if $c$ is a continuous constant, then $H^{+}(x)=H^-(x)$ by \eqref{interface1},
and thus the immersed interface cubic polynomial $H(x)$ is identical to the cubic interpolated polynomial \eqref{cubic}.

Let us consider the advection equation $u_t - (cu)_x = 0$, with $c>0$.
Assume that an interface locates $x=\alpha \in [x_{j-1},x_j]$. Then $x_{j-1}$
is an irregular point and the numerical solutions at $x_{j-1}$ are give by
\begin{equation*}
  u^{n+1}_{j-1} =H^-(x_{j-1}+c^-\Delta t),\quad v^{n+1}_{j-1}= H^-_x(x_{j-1}+c^-\Delta t).
\end{equation*}

\noindent
\textsf{Interface condition} $[c u]=0$.\\
The interface condition yields to the interface relations $[c^{i}\frac{\partial^iu}{\partial x^i}]=0$, $i\in\N$.
We define the immersed interface polynomial $H^{\pm}$ using the conditions
\begin{equation*}
[cu]=0,\quad [c^2u_x]=0,\quad [c^3u_{xx}]=0,\quad [c^4u_{xxx} ]=0,
\end{equation*}
and the interpolation condition \eqref{interp1}. The numerical solutions are given by the same form as
\begin{equation*}
  u^{n+1}_j =H^+(x_j-c^+\Delta t),\quad v^{n+1}_j= H^+_x(x_j-c^+\Delta t).
\end{equation*}
We list up the features and the advantages of IIM-CIP:
\begin{description}
\item [(i)] The method becomes the standard CIP if the discontinuities in
the coefficients disappear.
\item [(ii)] The structure of the IIM-CIP is as simple as CIP; The cubic
interpolation profile $H$ is replaced by the immersed interface cubic
polynomial.
\item [(iii)] The immersed interface cubic polynomial is an interpolation for the solution on the irregular cell, and the order of accuracy is four for $u(t,x)$ and three for its derivative.
    Therefore, the order of accuracy is maintained to be the same as that of the standard CIP scheme for $CFL \le 1$.
\item [(iv)] No grid refinement is required to maintain the accuracy in a vicinity of the interface.
%No special modification to the CIP scheme is needed to maintain of the order of accuracy.
\end{description}
\subsection{Numerical results}\label{subsec:numericalDISCON}
In this section we present some numerical results to illustrate IIM-CIP for discontinuous velocity and verify third order accuracy in time and space.\\
{\bf Example \ref{subsec:numericalDISCON}.1.} We consider the transport equation \eqref{equ:discon_tra} with discontinuity
\begin{equation*}
  c(x) = \left\{\begin{array}{ll}
           c_1 & 0   \le x \le \alpha, \\
           c_2 & \alpha \le x \le 1.
         \end{array}\right.
\end{equation*}
The periodic boundary condition is imposed. The interface condition $[u]=0$ and $[cu]=0$ at the interface $x=\alpha$ are tested. We apply IIM-CIP method developed in section \ref{sec:IIMCIP} to this problem. The constants $\alpha$, $c_1$ and $c_2$ are given below.
We take $u(0,x) = \exp(-\frac{(x-0.2)^2}{0.05^2})$, as initial condition. The exact solution for $[u]=0$ is given by $u(t,x) = u(0, y)$, where
\begin{equation*}
  y=\left\{
  \begin{array}{ll}
     x-c_2t, & \mbox{ if } x\ge \alpha \mbox{ and } t\le \frac{x-\alpha}{c_2},\\
     x-c_1t, & \mbox{ if } x \le \alpha, \\
     \frac{c_1}{c_2}(x-c_2 t) + (1-\frac{c_1}{c_2})\alpha, & \mbox{ otherwise }.
   \end{array}\right.
\end{equation*}
The exact solution for $[cu]=0$ is $u(t,x) = \frac{c(y)}{c(x)}u(0,y).$

In the numerical experiments, the spatial domain $[0,1]$ is uniformly discretized with mesh size $\Delta x = \frac{1}{N}$ for $N\in\{50, 100, 200, 400, 800, 1600\}$, and the time step size is $\Delta t = 0.5\Delta x$. The location of the interface is set to be $\alpha=0.5$. The velocity is set to be $c_1=1$ and $c_2=2$. The initial condition for the numerical solution $\{u^0_k\}_k^N$ is given by
$u^0_k=u(0,x_k)$, and $\{v^0_k\}_k^N$ is computed by central finite difference of $\{u^0_k\}_k^N$.
For each mesh size, the error in the numerical solutions at time $t=0.4$ is measured by \eqref{numericalerror}.

\begin{figure}[h]
\centering
\includegraphics[width=4cm, bb=124 252 500 540]{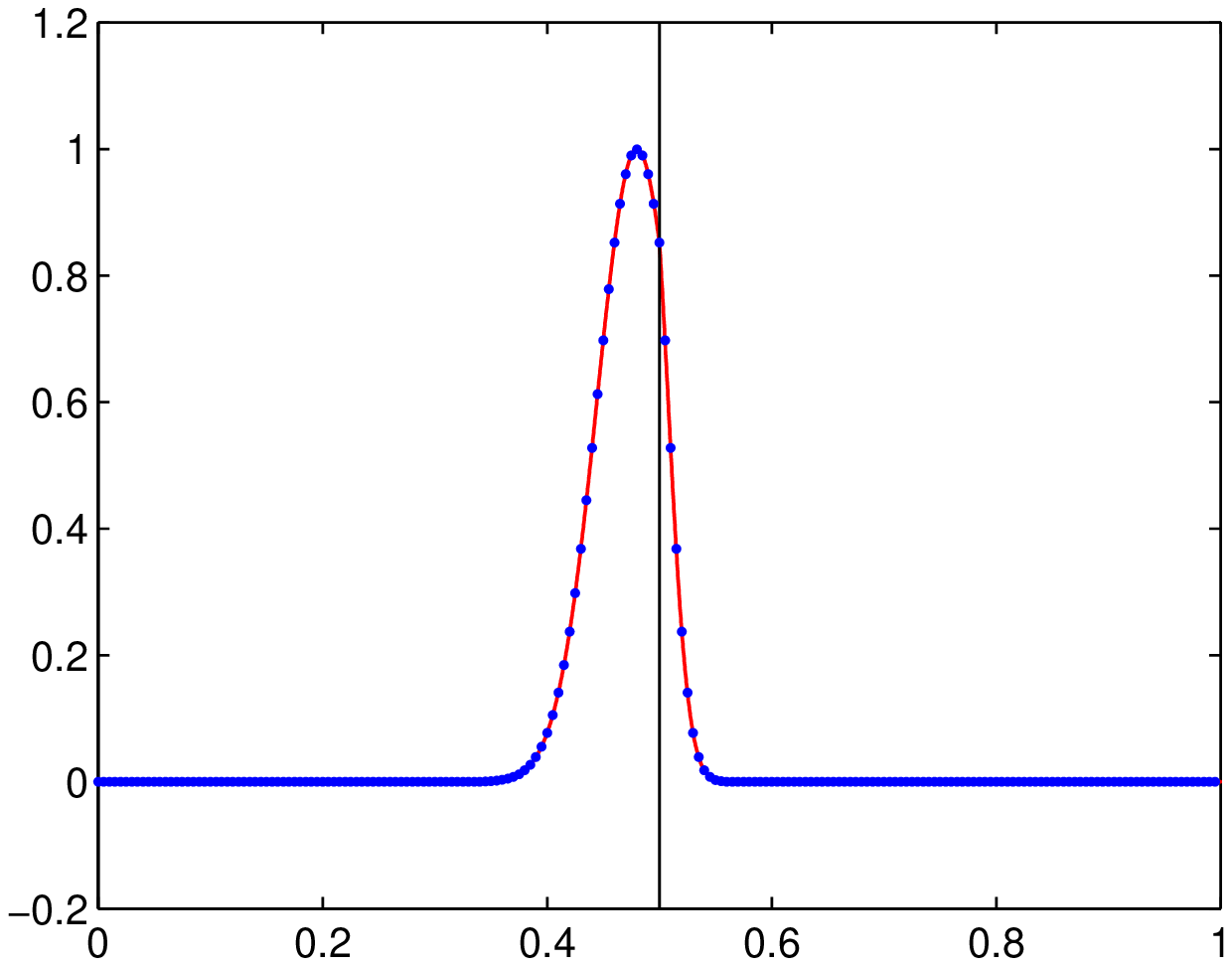}
\includegraphics[width=4cm, bb=124 252 500 540]{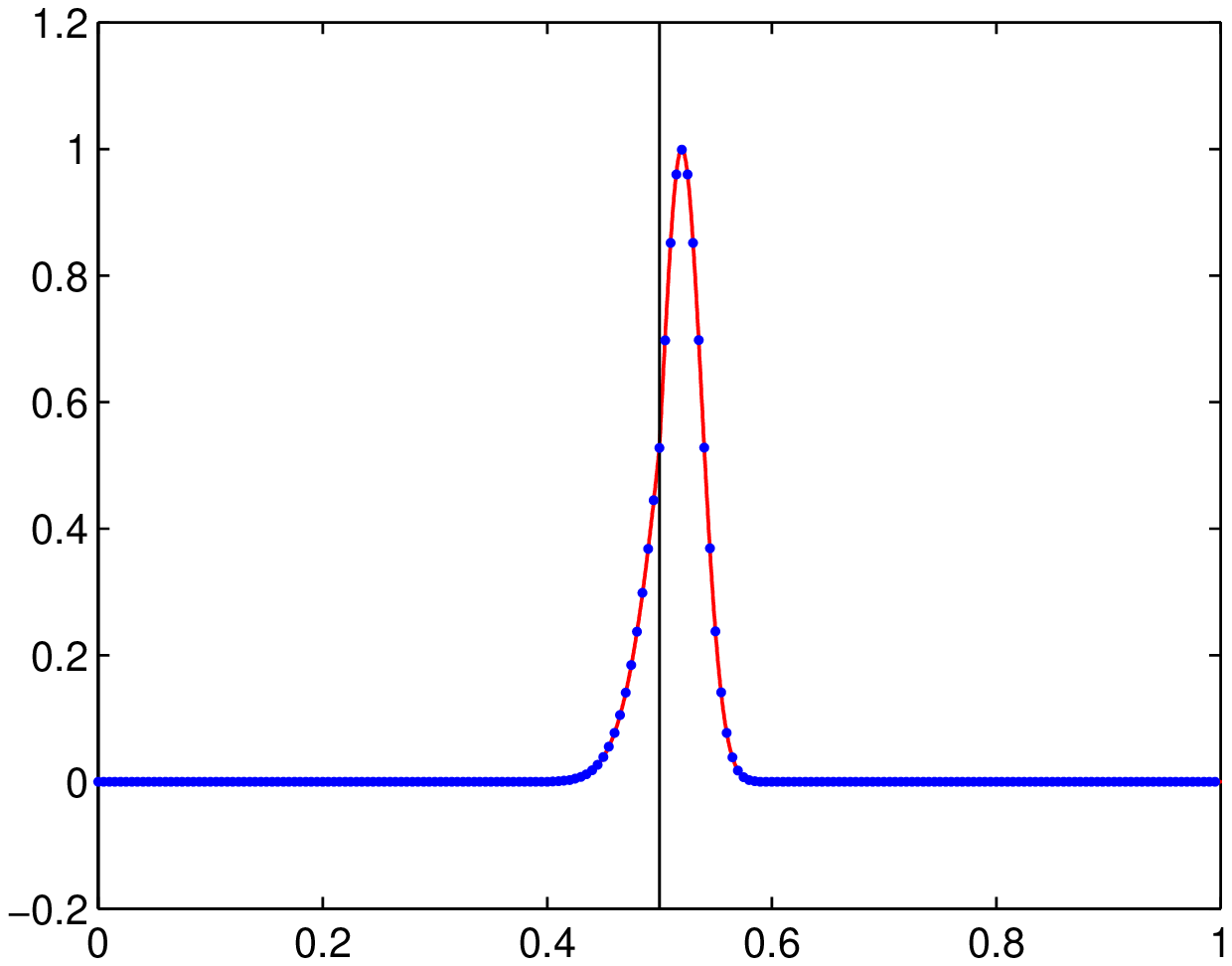}
\includegraphics[width=4cm, bb=124 252 500 540]{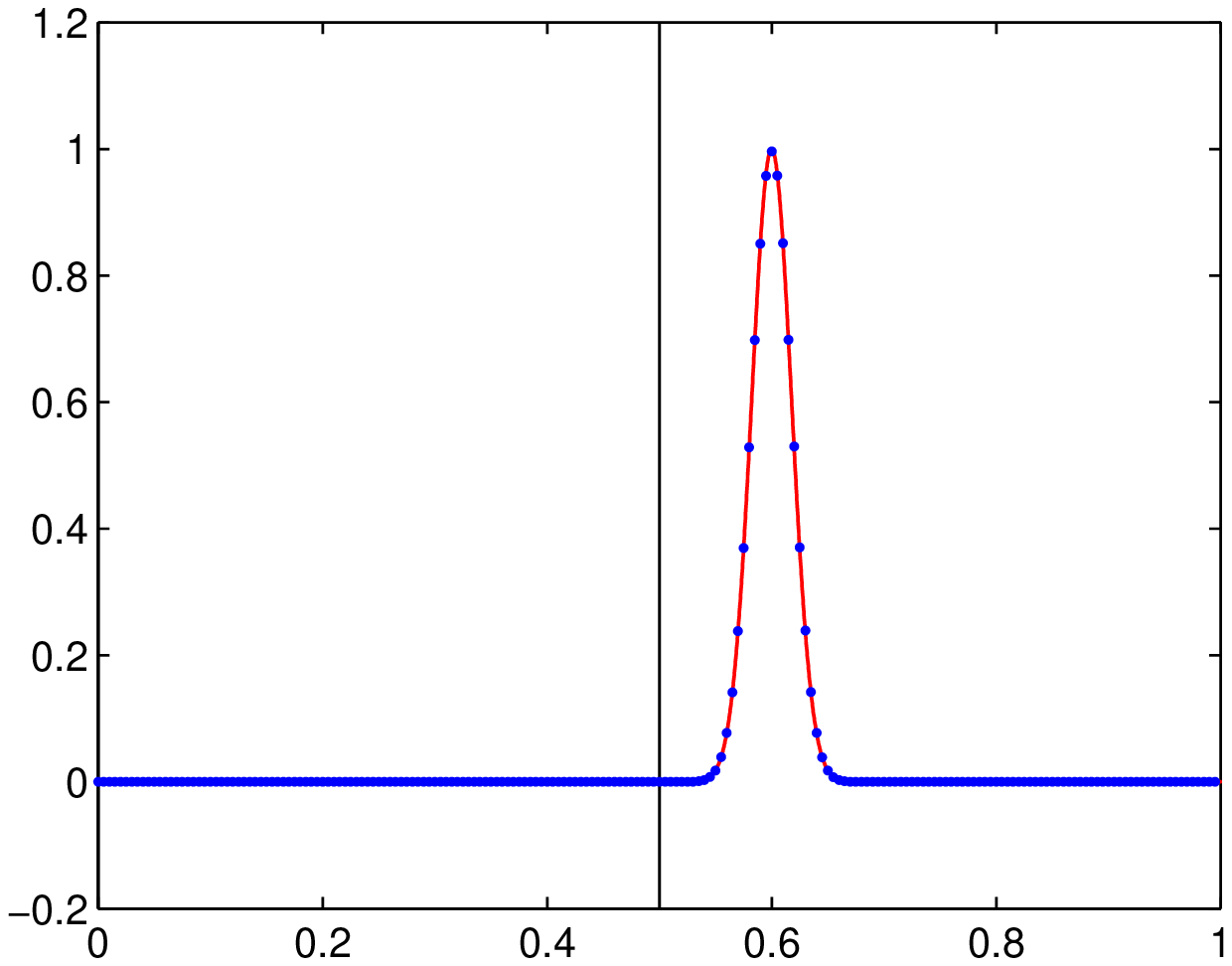}
\caption{1-D transport with $c^-=2$ and $c^+=1$ and jump condition $[u]=0$. The vertical solid line indicates the interface $\alpha=0.5$. The plots are the numerical solution (dot) and the exact solution (solid) at $t=0.14$ (left), $t=0.17$ (center), $t=0.25$ (right). The mesh size is $\frac{1}{200}$.}\label{fig:IIMCIPu}
\end{figure}

Plots of Figure \ref{fig:IIMCIPu} show numerical solutions (dot) and exact solutions (solid line) at time $t=0.14$ (left), $t=0.17$ (center) and $t=0.25$ (right) to the transport equation with jump condition $[u]=0$ at $\alpha=0.5$.
The mesh size $\frac{1}{200}$ is used to compute the numerical solutions.
The vertical solid line indicates the location of the interface. As the wave passes the interface, it slows down and becomes narrower. No spurious oscillation is observed in the numerical solutions.

\begin{figure}[h]
\centering
\includegraphics[width=4cm, bb=124 240 500 540]{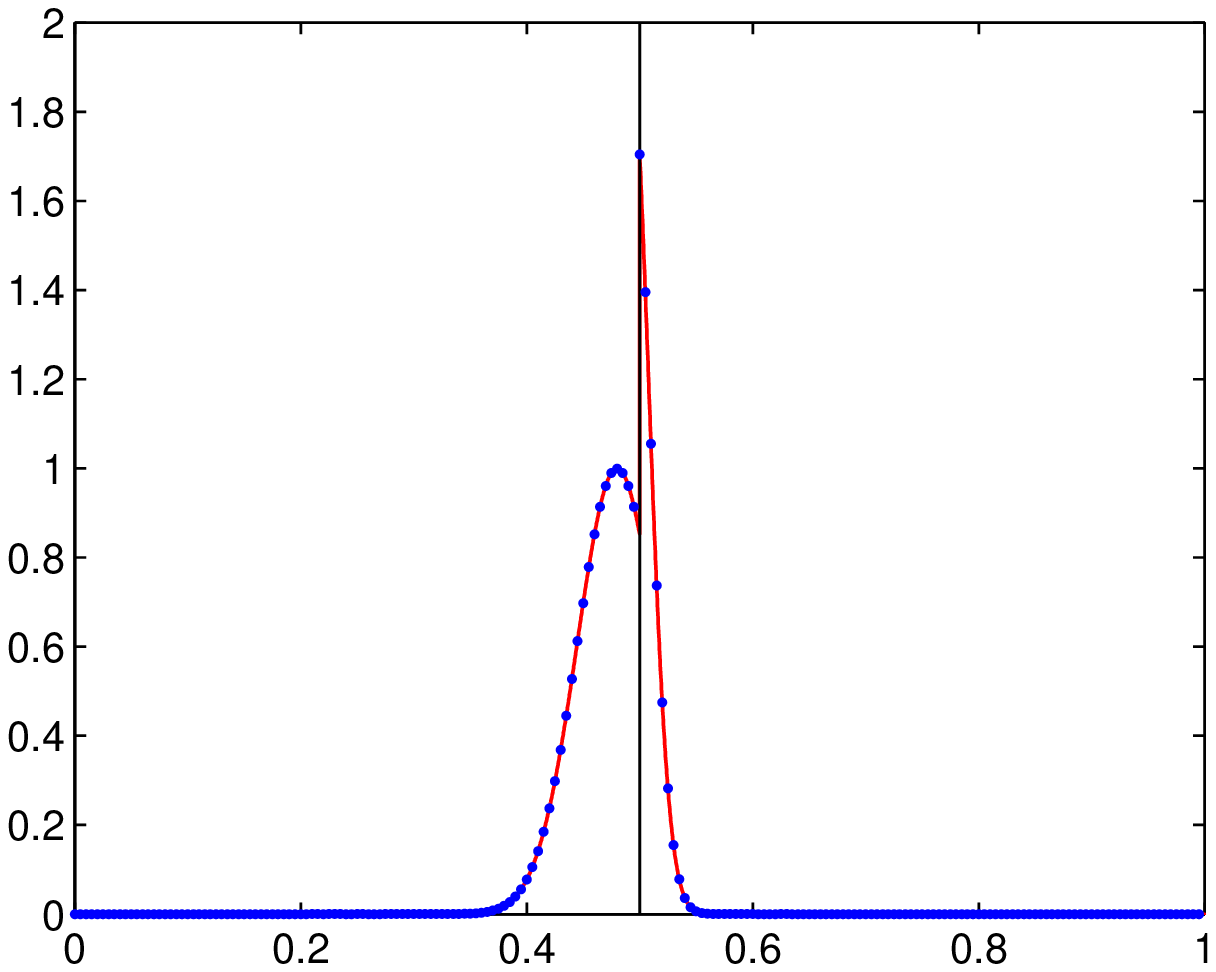}
\includegraphics[width=4cm, bb=124 240 500 540]{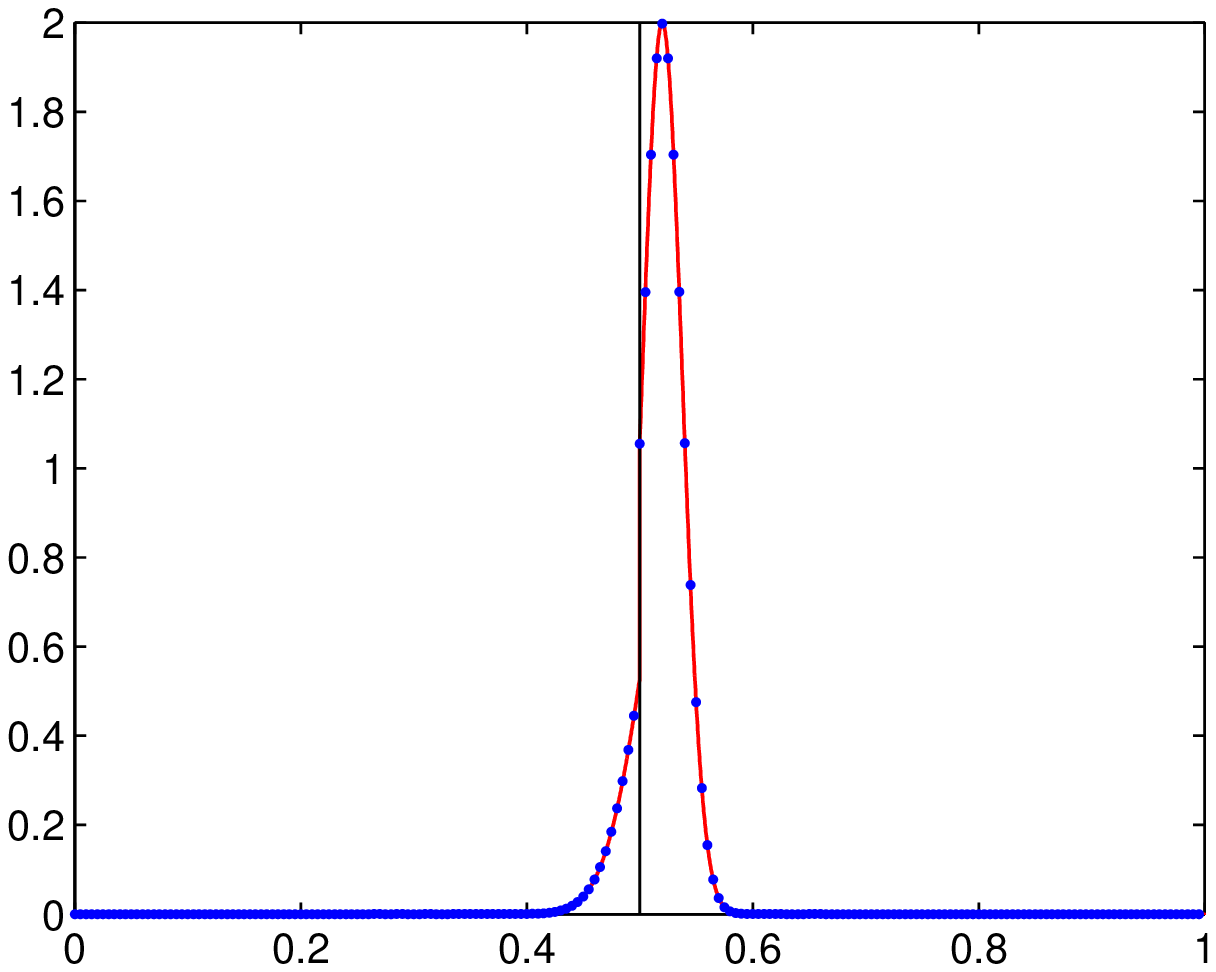}
\includegraphics[width=4cm, bb=124 240 500 540]{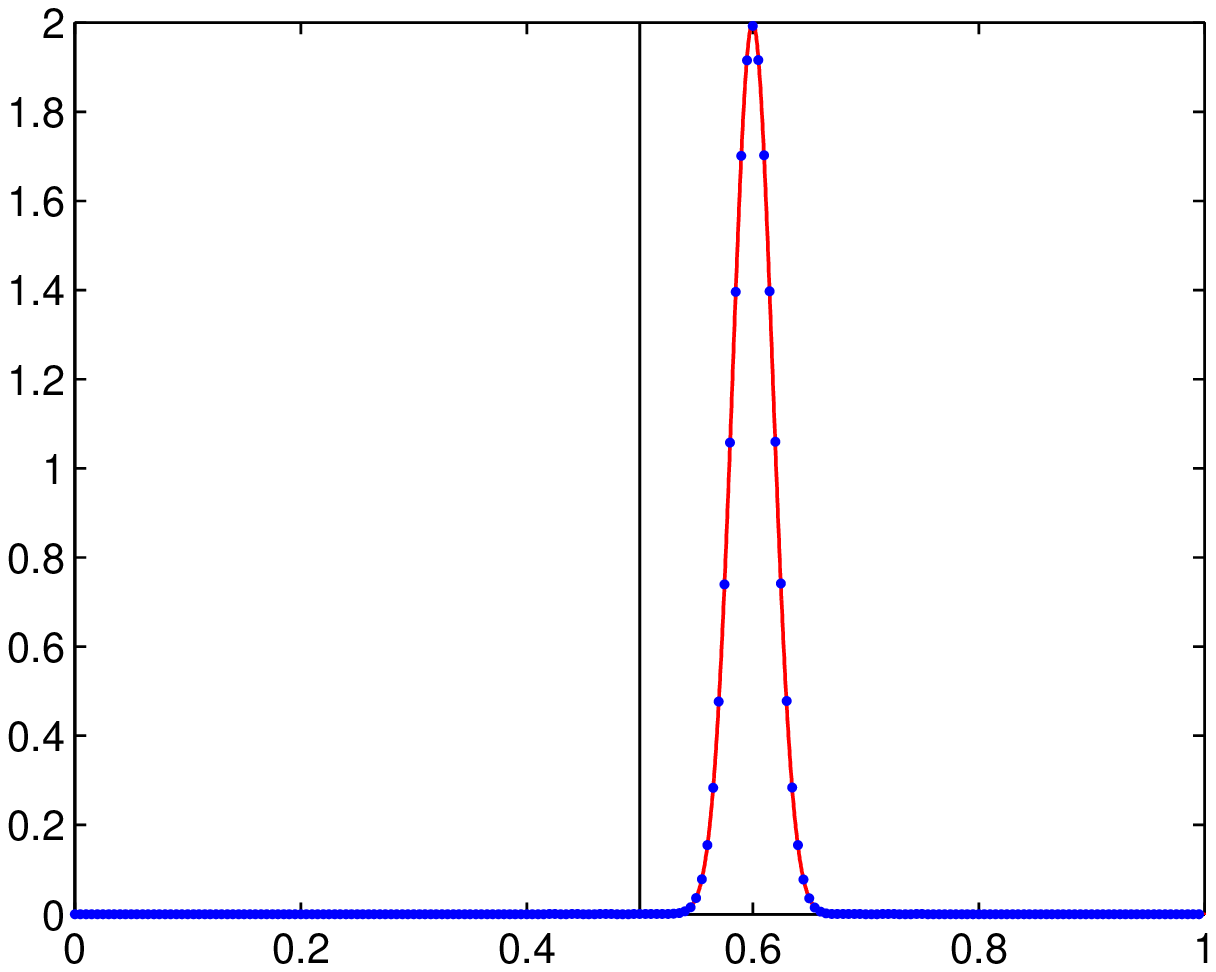}
\caption{1-D transport with $c^-=2$ and $c^+=1$ and jump condition $[cu]=0$. The vertical solid line indicates the interface $\alpha=0.5$. The plots are the numerical solution (dot) and the exact solution (solid) at $t=0.14$ (left), $t=0.17$ (center), $t=0.25$ (right). The mesh size is $\frac{1}{200}$.}\label{fig:IIMCIPcu}
\end{figure}

\begin{figure}[h]
\centering
\includegraphics[width=6cm,bb=124 252 500 540]{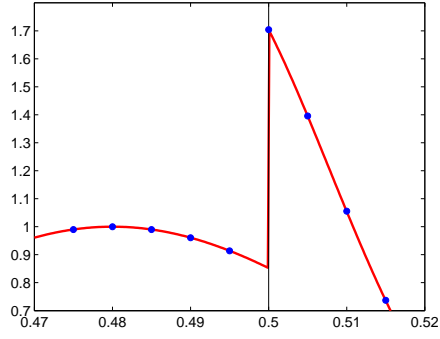}
\caption{A magnification of the left plot in figure \ref{fig:IIMCIPcu}.
The plots are the numerical solution (dot) and the exact solution (solid) at $t=0.14$.}\label{fig:IIMCIPcubig}
\end{figure}

Numerical solutions and the exact solution for the interface condition $[cu]$ is plotted in Figure \ref{fig:IIMCIPcu}.
The exact solution $u$ has a jump discontinuity at the interface.
We see that the numerical solution also exhibits the distinct jump at the interface.
A magnification of the figure at time $t=0.14$ around the interface is depicted in Figure \ref{fig:IIMCIPcubig} so that the jump discontinuity in the numerical solution and the exact solution are more visible.
We observe that the numerical solution almost coincides with the exact solution.

\begin{figure}[h]
\centering
\includegraphics[width=6cm, bb=100 230 500 540]{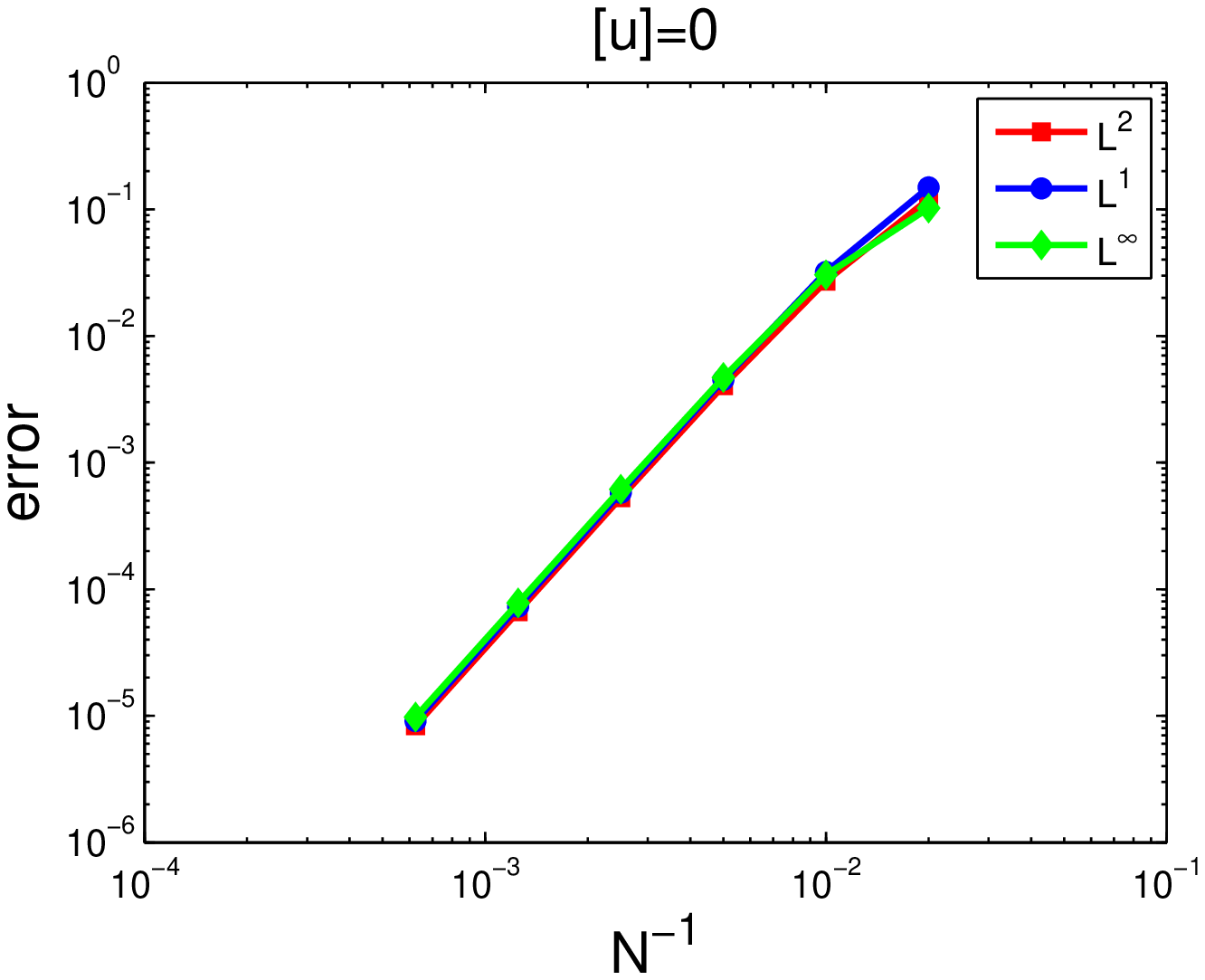}
\includegraphics[width=6cm, bb=100  230 500 540]{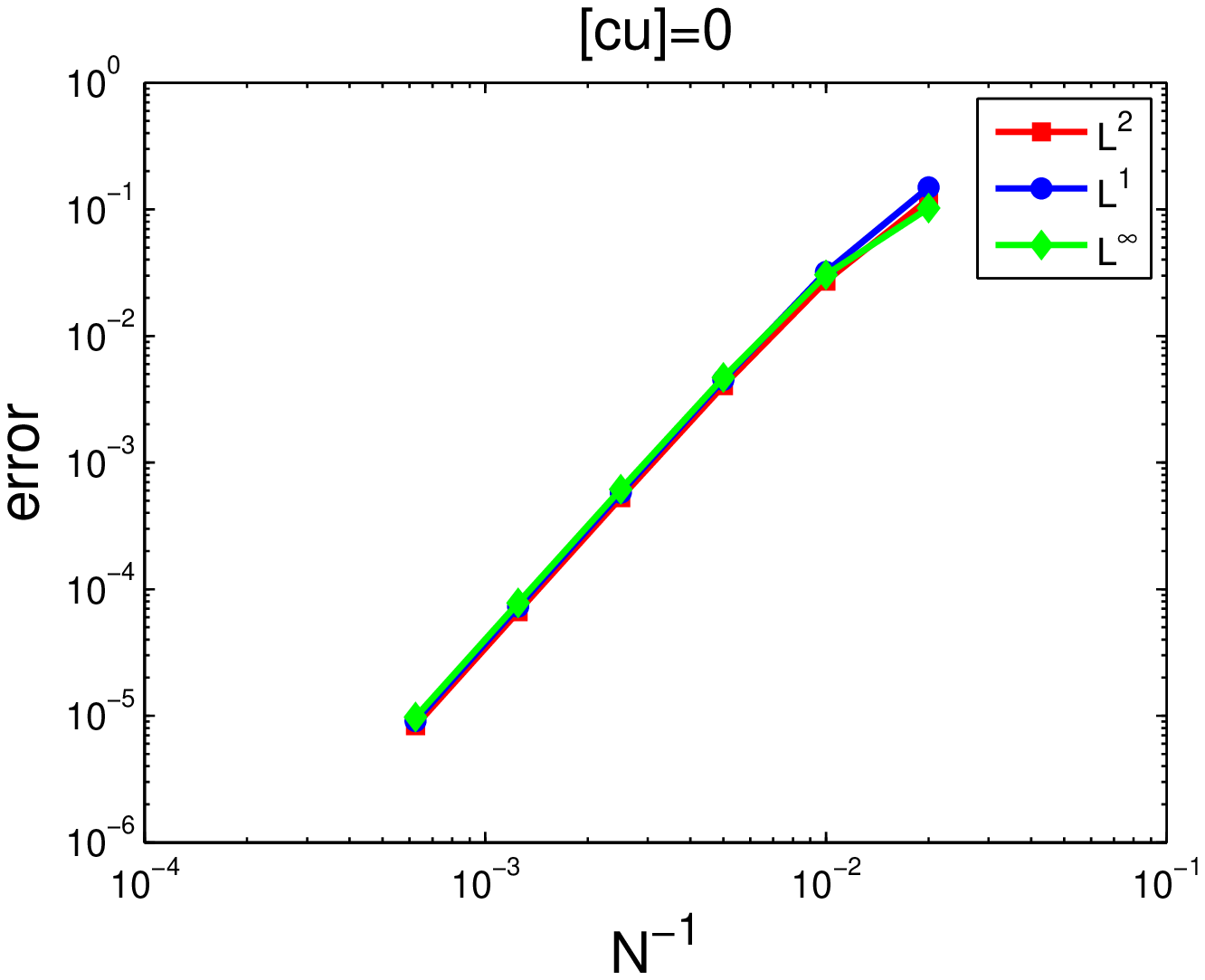}
\caption{The third order convergence against mesh size $N^{-1}$ of the numerical error in the numerical solution by the IIM-CIP. The error is computed by \eqref{numericalerror} at $t=0.4$.}
\label{fig:errorIIMCIPtraadv}
\end{figure}

Plots of Figure \ref{fig:errorIIMCIPtraadv} show the error in the computed solutions at $t=0.4$ against mesh size $N^{-1}$. Grid refinement studies confirm that the third order accuracy in time and space is achieved at all grid points. The third-order accuracy in time and space of the standard CIP is maintained even in the presence of the interface.
%This is explanted by the fact that the CIP scheme has small phase errors and small gain errors over a wide frequency range \cite{Utsumi+KunugiETAL-StabaccuCubiInte:97}.

In \cite{Zhang+LeVeque-immeintemethacou:97}, an immersed interface method is presented for a piecewise constant velocity. A piecewise quadratic interpolation on a cell $[x_{j-1},x_j]$
which contains a jump discontinuity in $c(x)$ is constructed based on the immersed interface method
and solution values at local three grid points. The underlined time integration method used in \cite{Zhang+LeVeque-immeintemethacou:97} is Lax-Wendorff method. The method causes visible oscillations in the numerical solution due to the fact that the Lax-Wendroff scheme is dispersive.

\section{IIM-CIP for Maxwell's equations in one dimension}\label{sec:IIMCIPmax}
Let us consider one dimensional Maxwell's equations
\begin{equation}\label{equ:max}
\begin{array}{ll}
\varepsilon E_t &= H_x\\
\mu H_t &= E_x
\end{array}
\end{equation}
for $x\in[0,1]$, $t>0$, with periodic boundary condition. %The coefficients $\varepsilon$ and $\mu$ are constant.

We begin with the case $\mu$ and $\varepsilon$ are constants. The solution update scheme is based on the D'Alambert formula:
\begin{equation*}
\begin{array}{l}
H(t_{n+1},x) =
    \tfrac{H(t_n,x-c\Delta t)+ H(t_n,x+c\Delta t)}{2}
-\tfrac{E(t_n,x-c\Delta t)- E(t_n,x+c\Delta t)}{2c\mu}, \\
E(t_{n+1},x) =
    \tfrac{E(t_n,x-c\Delta t)+ E(t_n,x+c\Delta t)}{2}
-\tfrac{H(t_n,x-c\Delta t)- H(t_n,x+c\Delta t)}{2c\epsilon},
\end{array}
\end{equation*}
where $c=\frac{1}{\sqrt{\mu\epsilon}}$. The formula follows from the fact that the transformations
\begin{equation*}
  u^1 = \sqrt{\mu}H - \sqrt{\varepsilon}E,\quad
  u^2 = \sqrt{\mu}H + \sqrt{\varepsilon}E,
\end{equation*}
satisfy the transport equations
\begin{equation*}
   u_t^1+cu_x^1= 0,\quad
   u_t^2+cu_x^2= 0.
\end{equation*}

 Let us denote the numerical solution to $H(t_n,x_k)$, $H_x(t_n,x_k)$ $E(t_n,x_k)$ and $E_x(t_n,x_k)$ by $H^{n}_{k}$, $DH^{n}_{k}$, $E^{n}_{k}$ and $DE^{n}_{k}$, respectively. Let $h^{k-1,k}(x)$ and $e^{k-1,k}(x)$ be the cubic polynomials on the grid $[x_{k-1},x_k]$ at time $t=t_n$ determined by the conditions
\begin{equation*}
   \begin{array}{llll}
   h^{k-1,k}(x_{k-1})=H^{n}_{k-1}, & h_x^{k-1,k}(x_{k-1}) =DH^{n}_{k-1} & h^{k-1,k}(x_{k})=H^{n}_{k}& h_x^{k-1,k}(x_{k}) = DH^{n}_{k}, \\
   e^{k-1,k}(x_{k-1})=E^{n}_{k-1}, & e_x^{k-1,k}(x_{k-1}) =DE^{n}_{k-1} & e^{k-1,k}(x_{k})=E^{n}_{k}& e_x^{k-1,k}(x_{k}) = DE^{n}_{k}.
   \end{array}
\end{equation*}
Let us assume the CFL number is less or equal to 1.
Then $x_k-c\Delta t \in [x_{k-1},x_k]$ and CIP scheme for equation \eqref{equ:max} is given by%The numerical solutions at next time level are then given by
\begin{equation*}
\begin{array}{ll}
H^{n+1}_k &= \tfrac{h^{k-1,k}(x_k-c\Delta t)+ h^{k,k+1}(x_k+c\Delta t)}{2}
-\tfrac{e^{k-1,k}(x_k-c\Delta t)-e^{k,k+1}(x_k+c\Delta t)}{2c\mu},\\
%DH^{n+1}_k &= \tfrac{h_x^{k-1,k}(x_k-c\Delta t) + h_x^{k,k+1}(x_k + c\Delta t)}{2}-\tfrac{e_x^{k-1,k}(x_k-c\Delta t)-e_x^{k,k+1}(x_k+c\Delta t)}{2c\mu},\\
E^{n+1}_k & = \tfrac{e^{k-1,k}(x_k-c\Delta t) + e^{k,k+1}(x_k + c\Delta t)}{2}
-\tfrac{h^{k-1,k}(x_k-c\Delta t)-h^{k,k+1}(x_k+c\Delta t)}{2c\varepsilon},%\\%\label{iimnsolE}
%DE^{n+1}_k &= \tfrac{e_x^{k-1,k}(x_k-c\Delta t) +e_x^{k,k+1}(x_k+c\Delta t)}{2}-\tfrac{h_x^{k-1,k}(x_k-c\Delta t)-h_x^{k,k+1}(x_k+c\Delta t)}{2c\varepsilon}.
\end{array}
\end{equation*}
and $DH^{n+1}_k$ and $DE^{n+1}_k$ are given by replacing the polynomials within $H^{n+1}_k$ and $E^{n+1}_k$ with the derivatives.

Now we consider the Maxwell's equations with $\varepsilon$ and $\mu$ being piecewise constants, i.e.,
$\varepsilon=\varepsilon^-$ and $\mu=\mu^-$ in $x<\alpha$ and $\varepsilon=\varepsilon^+$ and
$\mu=\mu^+$ in $\alpha<x$. The interface conditions for $H$ and $E$ at the interface are
$[H]=0$ and $[E]=0$. The interface relations are
\begin{eqnarray*}
\left[H\right]&=&0, \quad
\left[\frac{1}{\varepsilon} H_x\right] =0, \quad
\left[\frac{1}{\mu\varepsilon}H_{xx}\!\right]=0, \quad \left[\frac{1}{\mu\varepsilon^2}H_{xxx}\!\right]=0, \label{interfaceH11} \\
\left[E\right]&=&0, \quad
\left[\frac{1}{\mu}E_x\right]=0, \quad
\left[\frac{1}{\mu\varepsilon}E_{xx} \right]=0, \quad
\left[\frac{1}{\mu^2\varepsilon}E_{xxx} \right]=0, \label{interfaceE11}
\end{eqnarray*}
which are obtained in a usual manner in IIM; first differentiating the relation $[H]=0$ with respect to $t$,
then substituting the equation $\mu H_t=E_x$ to get $\left[\frac{1}{\mu}E_x\right] =0$.
Differentiate this relation again with respect to $t$ and substitute $\varepsilon E_t=H_x$,
we obtain $\left[\frac{1}{\mu\varepsilon} H_x\right] =0$. The others are given by repeating this procedure.

Let us assume that the interface $x=\alpha$ is included in
$[x_{j-1},x_j]$. %We update $H^{n}_j$ and $DH^n_j$ at an irregular
%point $x_{j}$ in the following way.
Let us denote the immersed interface cubic polynomials for
$H$ and $E$ on the cell $[x_{j-1},x_j]$ by $h^{\pm}(x)$ and
$e^{\pm}(x)$ respectively. From the interface relations,
\begin{align*}
  h^{\pm}(x)& = a_0 + \varepsilon^{\pm}a_1\frac{x-\alpha}{\Delta x} + (\mu\varepsilon)^{\pm}\frac{a_2}{2}(\frac{x-\alpha}{\Delta x})^2 +
  (\mu\varepsilon^2)^{\pm}\frac{a_3}{3!}(\frac{x-\alpha}{\Delta x})^3,\\
  e^{\pm}(x)& = b_0 + \mu^{\pm}b_1\frac{x-\alpha}{\Delta x} + (\mu\varepsilon)^{\pm}\frac{b_2}{2}(\frac{x-\alpha}{\Delta x})^2 +
  (\mu^2\varepsilon)^{\pm}\frac{b_3}{3!}(\frac{x-\alpha}{\Delta x})^3.
\end{align*}
The coefficients $a$ and $b$ are determined by the interpolation condition at two end points $x_{j-1},\; x_j$, i.e.,
\begin{equation*}
  \begin{array}{l}
A(\epsilon,\mu)a =(H^n_{j},\Delta x\; DH^n_{j}, H^n_{j-1},\Delta x\; DH^n_{j-1})^\top,\\
A(\mu,\epsilon)b =(E^n_{j},\Delta x\; DE^n_{j}, E^n_{j-1},\Delta x\; DE^n_{j-1})^\top
  \end{array}
\end{equation*}
%i.e., $ A(\epsilon,\mu)a=f$, $ A(\mu,\epsilon)b = g$
where
\begin{align*}
  A(\epsilon,\mu)=\left(
      \begin{array}{cccc}
        1 & \varepsilon^{+}{\theta} & \frac{(\mu\varepsilon)^{+}\theta^2}{2} & \frac{(\mu\varepsilon^2)^{+}\theta^3}{3!} \\
        0 & \varepsilon^{+} & (\mu\varepsilon)^{+}{\theta} & \frac{(\mu\varepsilon^{2})^+\theta^2}{2} \\
        1 & \varepsilon^{-}(\theta-1) & \frac{(\mu\varepsilon)^{-}(\theta-1)^2}{2} & \frac{(\mu\varepsilon^2)^{-}(\theta-1)^3}{3!} \\
        0 & \varepsilon^{-} & (\mu\varepsilon)^{-}(\theta-1)
        & \frac{(\mu\varepsilon^{2})^-(\theta-1)^2}{2}
      \end{array}
    \right),\quad \theta=\frac{x_{j}-\alpha}{\Delta x}.
%    f=\left(
%        \begin{array}{c}
%          H^n_{j} \\
%         \Delta x \;DH^n_{j} \\
%          H^n_{j-1} \\
%          \Delta x\; DH^n_{j-1}
%        \end{array}
%      \right),
\end{align*}

The numerical solutions $H^{n+1}_j$ and
$E^{n+1}_j$ at the irregular point $x_j$ are then given by
\begin{equation}\label{EMIIM}
\begin{array}{ll}
  H^{n+1}_j  &= \tfrac{h^+(x_j-c^+\Delta t)+ h^{j,j+1}(x_j+c^+\Delta t)}{2}
-\tfrac{e^+(x_j-c^+\Delta t)-e^{j,j+1}(x_j+c^+\Delta t)}{2c^+\mu^+},\\
%DH^{n+1}_j &= \tfrac{h_x^{+}(x_j-c^+\Delta t) + h_x^{j,j+1}(x_j + c^+\Delta t)}{2}-\tfrac{e^+_x(x_j-c^+\Delta t)-e_x^{j,j+1}(x_j+c^+\Delta t)}{2c^+\mu^+},\\
E^{n+1}_j &= \tfrac{e^+(x_j-c^+\Delta t) + e^{j,j+1}(x_j + c^+\Delta t)}{2}
-\tfrac{h^+(x_j-c^+\Delta t)-h^{j,j+1}(x_j+c^+\Delta t)}{2c^+\varepsilon^+},%\\ %\label{iimnsolE}
%DE^{n+1}_j &= \tfrac{e_x^+(x_j-c^+\Delta t) +e_x^{j,j+1}(x_j+c^+\Delta t)}{2}-\tfrac{h^+_x(x_j-c^+\Delta t)-h_x^{j,j+1}(x_j+c^+\Delta t)}{2c^+\varepsilon^+}.%\label{iimnsolDE}
\end{array}
\end{equation}
where $c^+=\frac{1}{\sqrt{\mu^+\epsilon^+}}$. $DH^{n+1}_j$ and $DE^{n+1}_j$ are given by replacing the polynomials within $H^{n+1}_j$ and $E^{n+1}_j$ with the derivatives.
We note that the point $x_{j-1}$ is also an irregular point, and the immersed interface cubic polynomials are used to update the numerical solutions at this point. For instance, $H_{j-1}^{n+1}$ is given as
 \begin{align*}
  H^{n+1}_{j-1} = \tfrac{h^{j-2,j-1}(x_{j-1}-c^-\Delta t)+ h^-(x_{j-1}+c^-\Delta t)}{2}
-\tfrac{e^{j-2,j-1}(x_{j-1}-c^-\Delta t)-e^-(x_{j-1}+c^-\Delta t)}{2c^-\mu^-},
\end{align*}
where $c^-=\frac{1}{\sqrt{\mu^-\epsilon^-}}$.
%
%One can refer to \cite{Barada+FukudaETAL-Cubiintepropsche:06} for another treatment of material discontinuity for CIP method.
\subsection{Application to Maxwell's equations with variable material parameters}\label{subsec:IIMCIPmax}
We apply the method we developed for variable
$\varepsilon(x),\;\mu(x)$. We approximate $\epsilon(x),\;\mu(x)$ by
the piecewise constant (discontinuous) media:
$$
\bar{\varepsilon}(x)=\tfrac{1}{\Delta
x}\int^{x_{j}+\frac{\Delta x}{2}}_{x_{j}-\frac{\Delta x}{2}}\varepsilon(x)\,dx,\quad
\bar{\mu}(x)=\tfrac{1}{\Delta
x}\int^{x_{j}+\frac{\Delta x}{2}}_{x_{j}-\frac{\Delta x}{2}}\mu(x)\,dx,
$$
on $(x_{j-1/2},x_{j+1/2})$ for all $j$. Then one can apply \eqref{EMIIM} for the both
backward and forward manner to obtain
\begin{equation}\label{EMIIMCIP}
\begin{array}{l}
H^{n+1}_j  = \tfrac{h^+_{j-1,j}(x_j-c\Delta t)+h^-_{j,j+1}(x_j+c\Delta t)}{2}
-\tfrac{e^+_{j-1,j}(x_j-c\Delta t)-e^-_{j,j+1}(x_j+c\Delta t)}{2c\mu},\\
E^{n+1}_j = \tfrac{e^+_{j-1,j}(x_j-c\Delta t) +
e^-_{j,j+1}(x_j + c\Delta t)}{2}
-\tfrac{h^+_{j-1,j}(x_j-c\Delta t)-h^-_{j,j+1}(x_j+c\Delta t)}{2c\varepsilon}.\\%\label{iimnsolE}
%\label{iimnsolDE}
%\label{iimnsolDH}\\
\end{array}
\end{equation}
By replacing the polynomials within \eqref{EMIIMCIP} with the derivatives, we obtain $DH^{n+1}_j$ and $DE^{n+1}_j$.
%DH^{n+1}_j = \tfrac{\partial_xh^+_{j-1,j}(x_j-c\Delta t) +
%\partial_xh^-_{j,j+1}(x_j + c\Delta t)}{2}
%-\tfrac{\partial_xe^+_{j-1,j}(x_j-c\Delta t)-\partial_xe^-_{j,j+1}(x_j+c\Delta t)}{2c\mu},\\
%DE^{n+1}_j = \tfrac{\partial_xe^+_{j-1,j}(x_j-c^+\Delta t)
%+\partial_xe^-_{j,j+1}(x_j+c^+\Delta t)}{2}
%-\tfrac{\partial_xh^+_{j-1,j}(x_j-c\Delta t)-\partial_xh^-_{j,j+1}(x_j+c\Delta t)}{2c\varepsilon}.
Here
$
c=\bar{c}_j,\; \mu=\bar{\mu}_j,\;\varepsilon=\bar{\varepsilon}_j
\quad \mbox{on the cell $(x_{j-1/2},x_{j+1/2})$},
$ and $h^\pm_{j-1,j}(x)$, $e^\pm_{j-1,j}(x)$ are the immersed interface cubic polynomials on $[x_{j-1},x_j]$ and
$h^\pm_{j,j+1}(x)$, $e^\pm_{j,j+1}(x)$ are the immersed interface cubic polynomials on $[x_{j},x_{j+1}]$.
%Here we have the two sided interface $x-\frac{\Delta x}{2}\in (x_{j-1},x_j)$ and
%$x_j+\frac{\Delta x}{2}\in (x_j,x_{j+1})$.

\subsection{Numerical results}\label{subsec:numericalEM}
We present two examples to illustrate the potential of the IIM-CIP for the Maxwell's equations. \\
{\bf Example \ref{subsec:numericalEM}.1.}
Consider the Maxwell's equations \eqref{equ:max} with
\begin{equation*}
 \varepsilon(x) = \left\{\begin{array}{ll}
           \varepsilon^- =1,& 0   \le x \le \alpha, \\
           \varepsilon^+ =\frac{4}{3},& \alpha \le x \le 1.
         \end{array}\right.,\quad
 \mu(x) = \left\{\begin{array}{ll}
           \mu^- =1,& 0   \le x \le \alpha, \\
           \mu^+ =3,& \alpha \le x \le 1.
         \end{array}\right.
\end{equation*}
The spatial domain $[0,1]$ is uniformly discretized with mesh size $\Delta x = \frac{1}{200}$, and the time step size is $\Delta t = 0.5\Delta x$. The location of the interface is set to be $\alpha=0.5$. As an initial condition, we take $H(0,x)=\exp(-\tfrac{(x-0.2)^2}{0.05^2})$ and $E(0,x)=-\sqrt{\frac{\mu}{\varepsilon}}H(0,x)$.

\begin{figure}[h]
\centering
\includegraphics[width=6cm, bb=124 240 500 570]{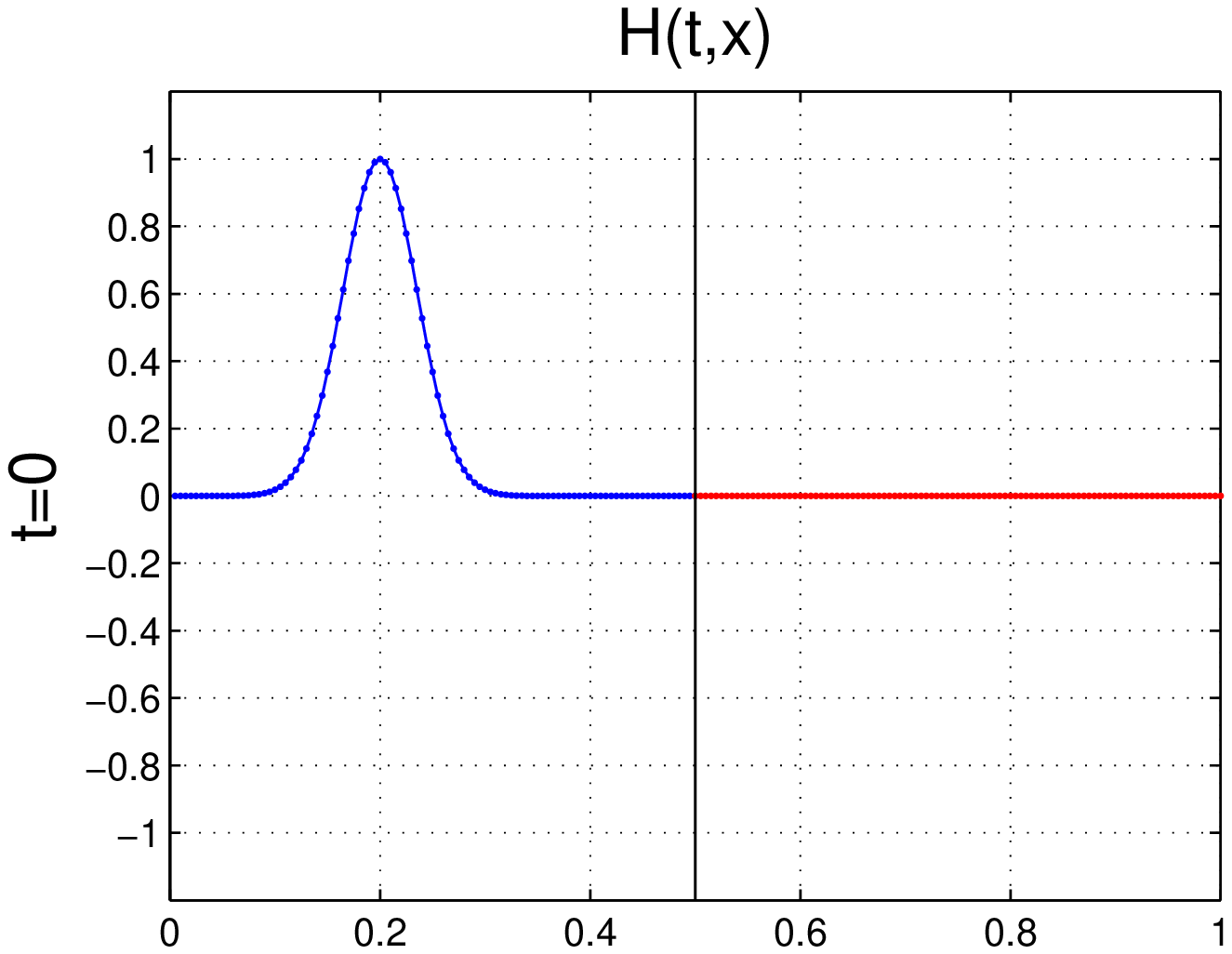}
\includegraphics[width=6cm, bb=124 240 500 570]{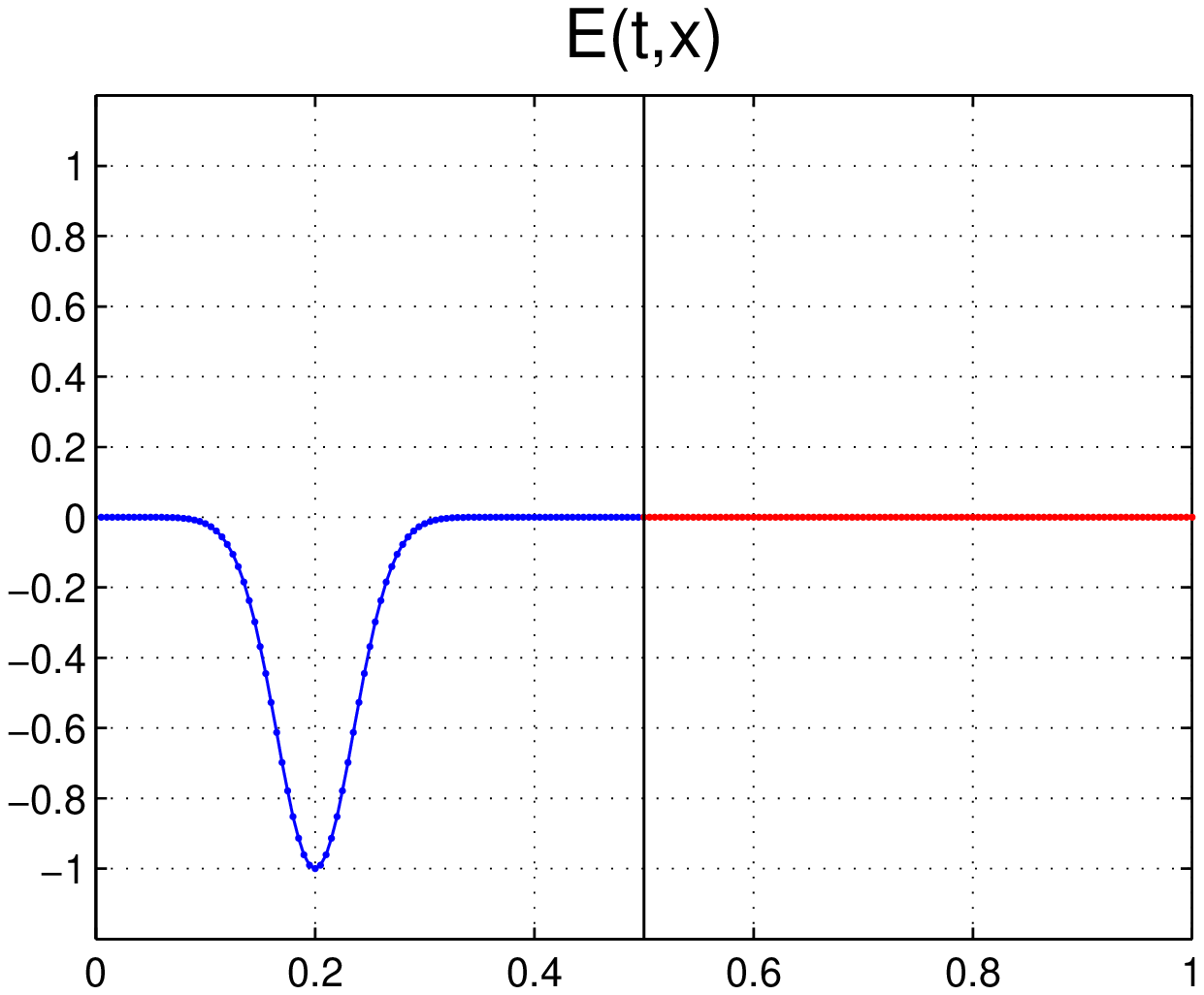}
\includegraphics[width=6cm, bb=124 240 500 540]{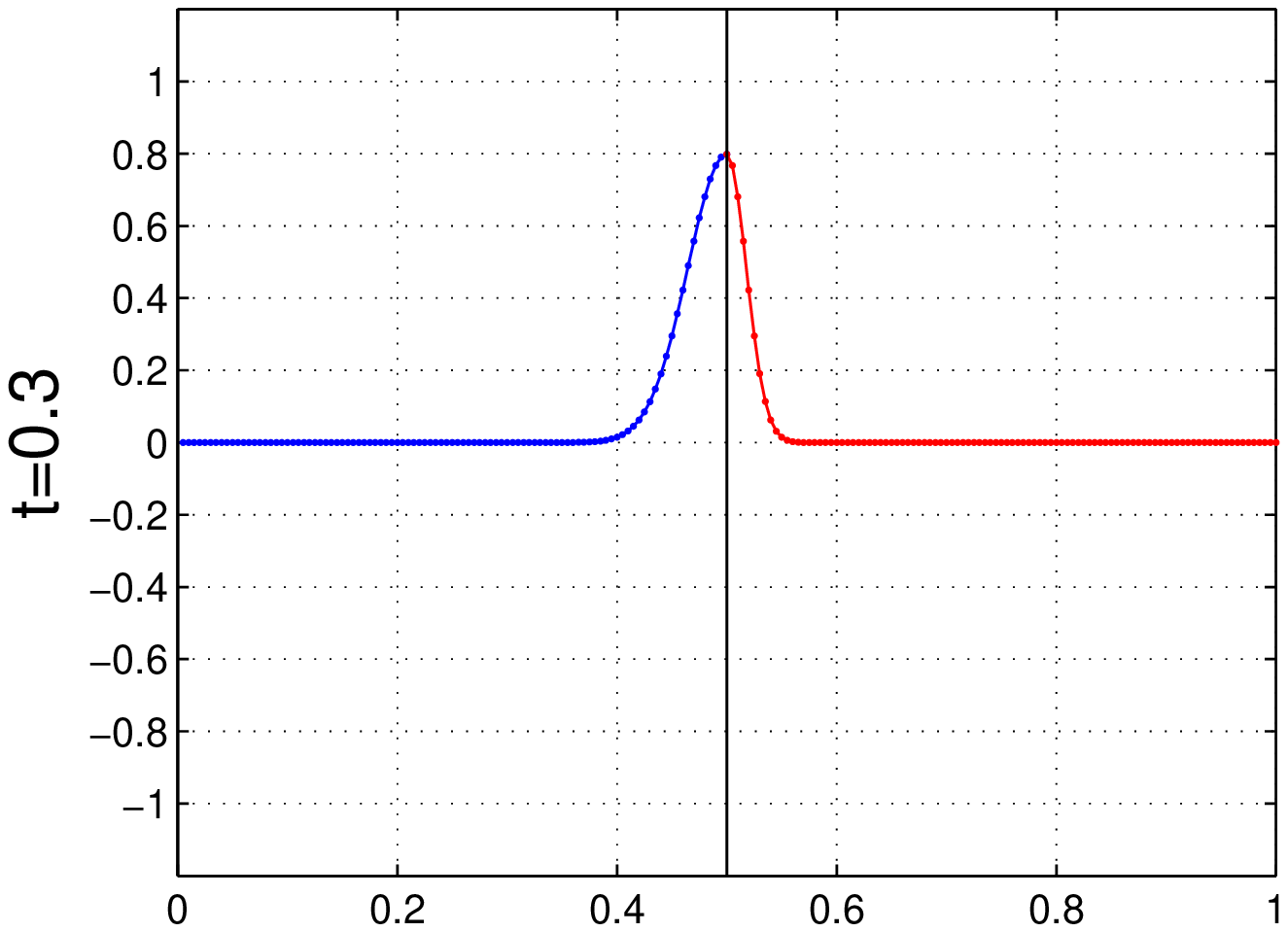}
\includegraphics[width=6cm, bb=124 240 500 540]{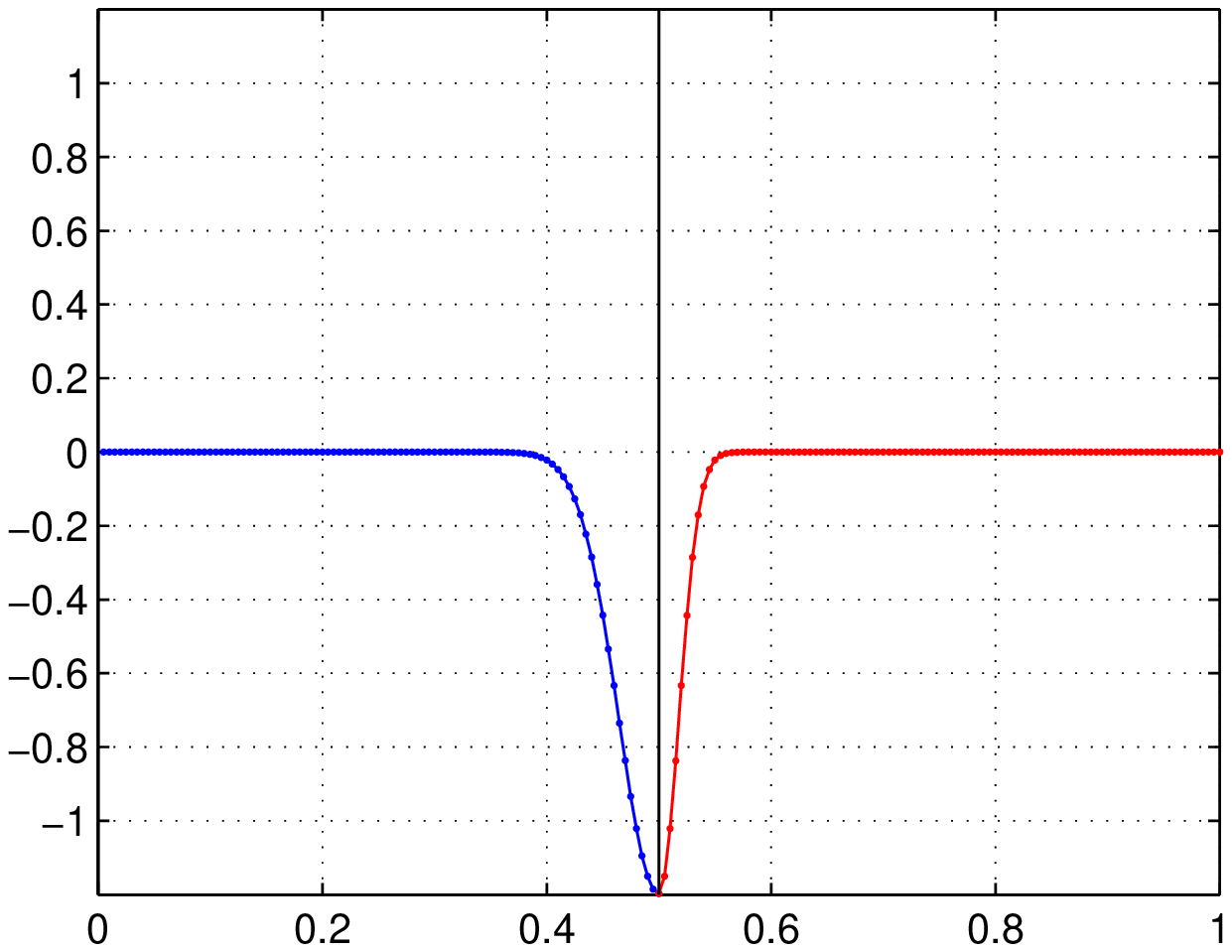}
\includegraphics[width=6cm, bb=124 240 500 540]{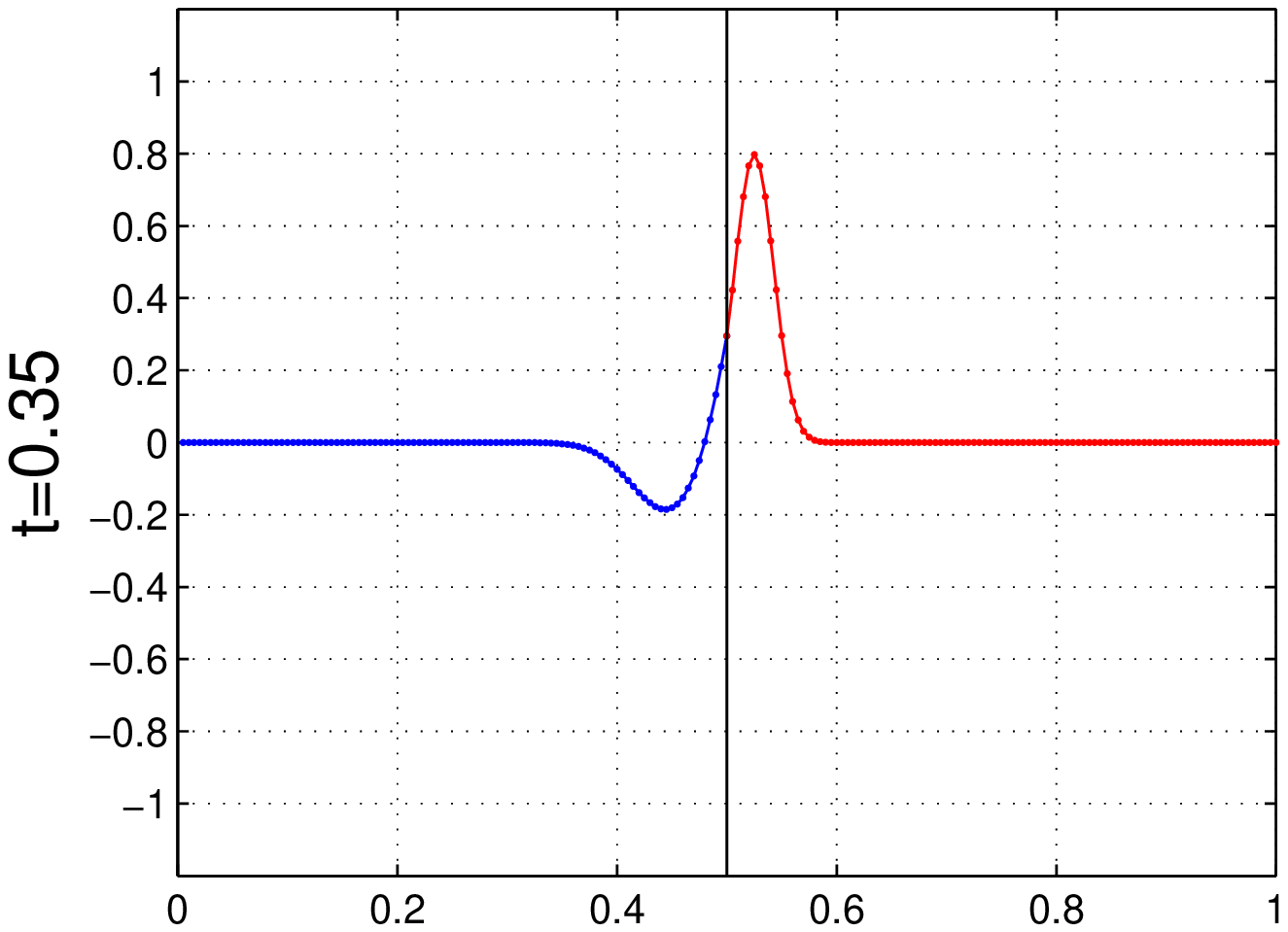}
\includegraphics[width=6cm, bb=124 240 500 540]{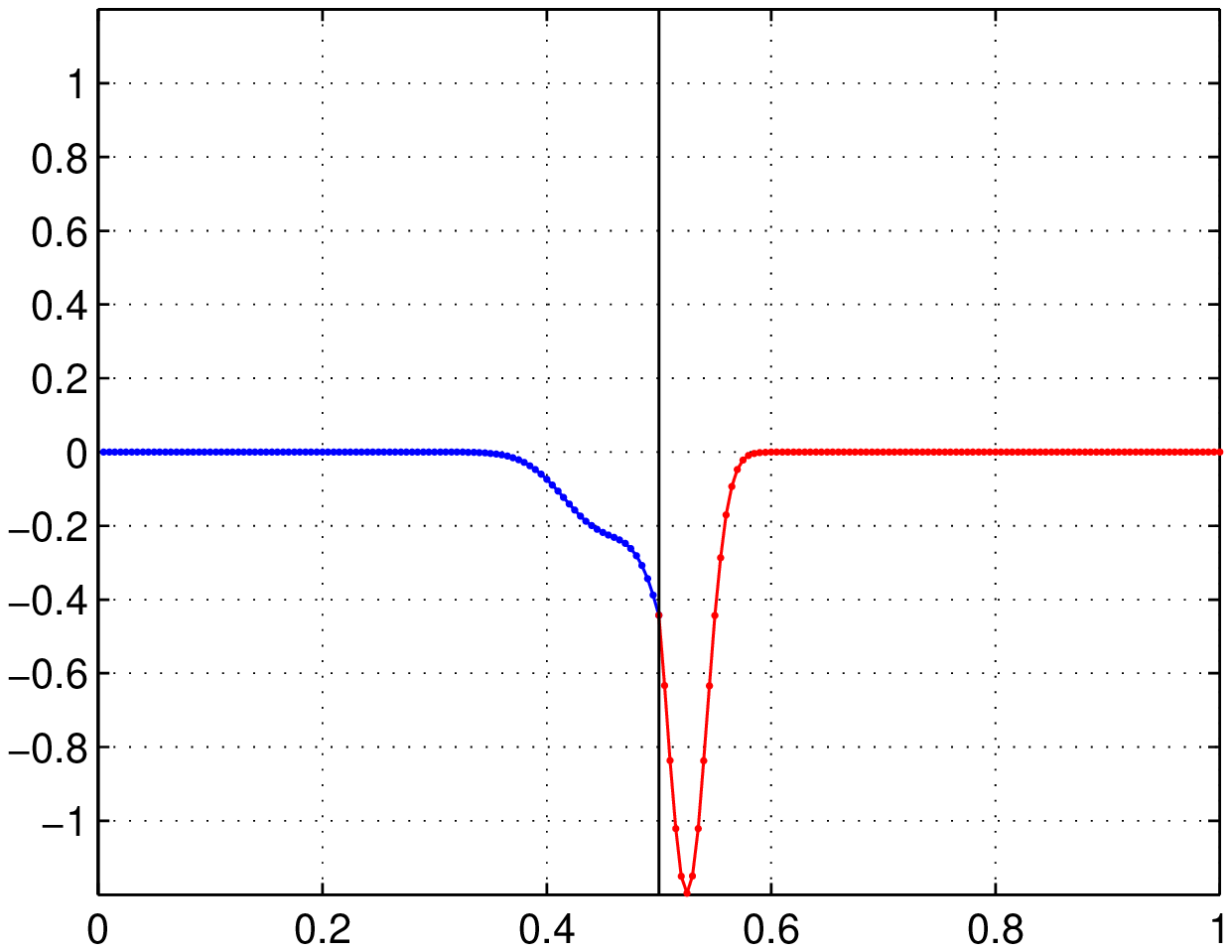}
\includegraphics[width=6cm, bb=124 240 500 540]{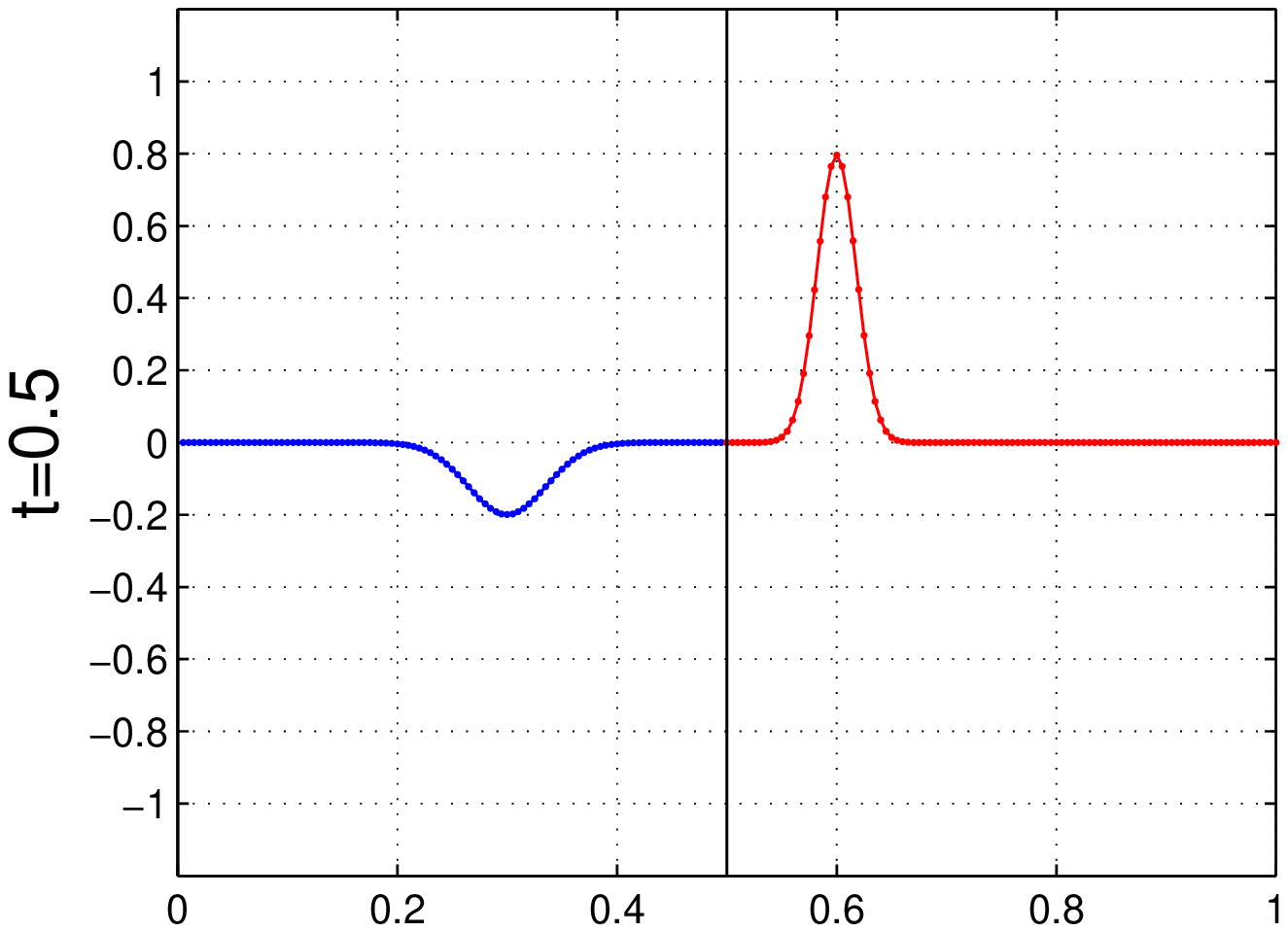}
\includegraphics[width=6cm, bb=124 240 500 540]{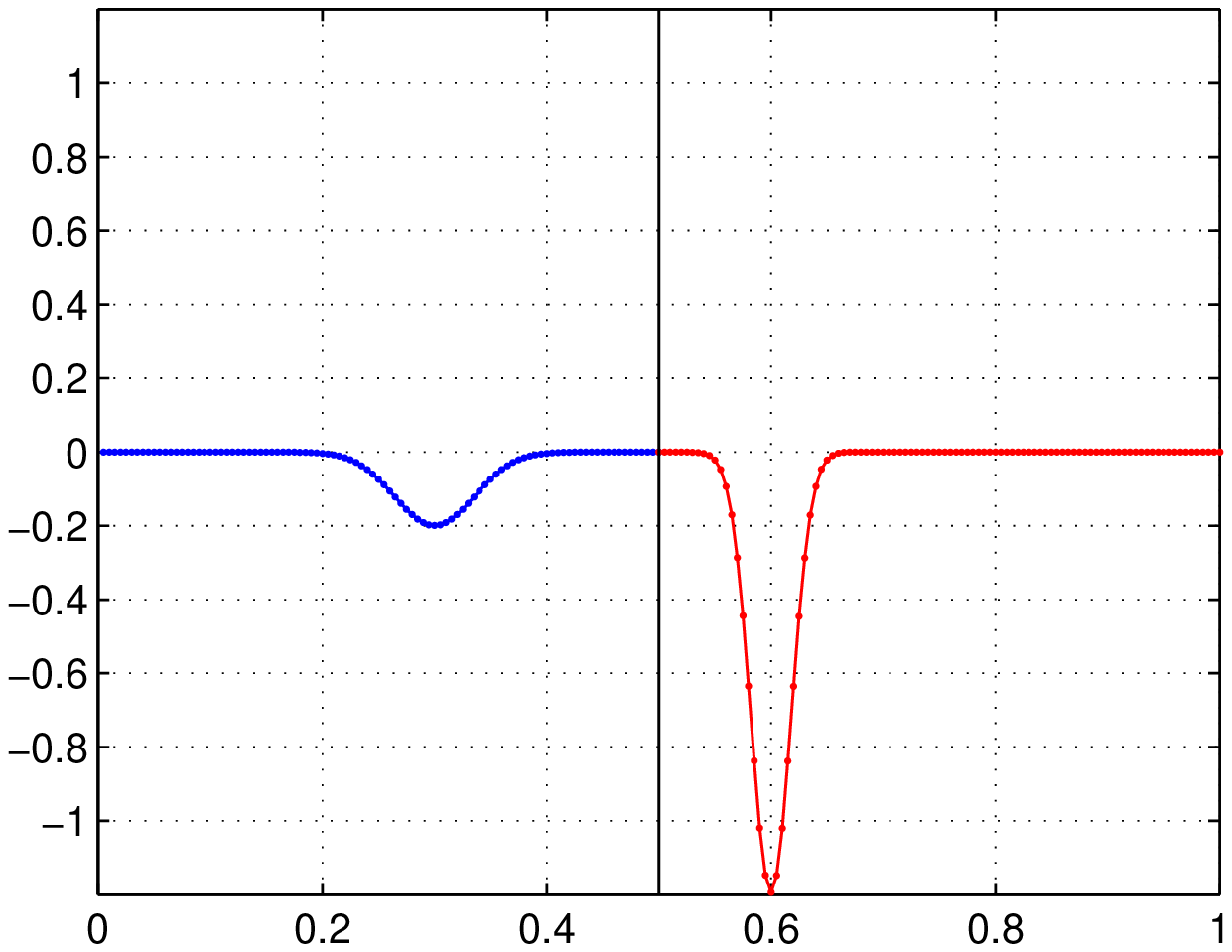}
\caption{1-D Maxwell's equations with $\mu^-=\varepsilon^-=1$ ($x<0.5$)and $\mu^+=3$, $\varepsilon^+=\frac{4}{3}$ ($0.5\le x$). The interface condition $[H]=[E]=0$ is imposed at the interface $x=0.5$ (vertical line). The left plots are snap shots of the numerical solution to $H$ at $t=0$, $t=0.3$, $t=0.35$ and $t=0.5$ from top to bottom. The right plots are the numerical solution to $E$ at the same time. The mesh size is $\frac{1}{200}$.}\label{fig:IIMCIPEM}
\end{figure}

Plots of Figure \ref{fig:IIMCIPEM} show the numerical solutions to $H(t,x)$ (left column) and $E(t,x)$ (right column)  at time $t=0$, $t=0.3$, $t=0.35$ and $t=0.5$.
There are no spurious oscillations observed in the vicinity of the interface, at least for this example.

\noindent
{\bf Example \ref{subsec:numericalEM}.2.}
Consider the Maxwell's equations
$\varepsilon(x) E_t = H_x$, $\mu(x) H_t = E_x$ for $x\in[0,1], \; t>0$ with periodic boundary condition, where $\varepsilon(x)=\mu(x)=\frac{1}{2}\cos(4\pi x) + 1$.
As an initial condition, we take $H(0,x)=\exp(-\tfrac{(x-0.5)^2}{0.05^2})$ and $E(0,x)=0$.
We apply the IIM-CIP developed in Section \ref{subsec:IIMCIPmax}. In this numerical test, the time step size is chosen to be $\Delta t= \frac{0.5 N^{-1}}{\max_{x\in[0,1]} c(x) } = 0.25  N^{-1}$ for each $N\in\{50, 100, 200, 400, 800, 1600\}$. The numerical solution is integrated in time by \eqref{EMIIMCIP}.

\begin{figure}[h]
\centering
\includegraphics[width=0.5\textwidth, bb=124 252 500 540]{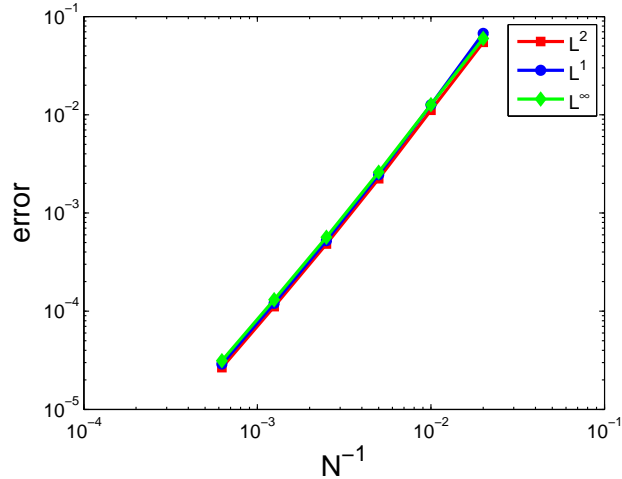}
\caption{The error \eqref{errorEM} in the numerical solutions at $t=1$ of the Maxwell's equation with variable material parameters $\varepsilon(x)=\mu(x) = \frac{1}{2}\cos(4\pi x) + 1$ computed by IIM-CIP against mesh size $N^{-1}$. The second order convergence in time and space is observed.}\label{fig:EMIIMCIP}
\end{figure}

The numerical solutions of the magnetic field $H$ at time $t=1$ are compared to the exact solution, which is identical to the initial condition, i.e., $H(1,x) = H(0,x)$. For each mesh size, the error in the numerical solutions is measured by $\ell^1$, $\ell^2$ and $\ell^\infty$ norm:
\begin{equation}\label{errorEM}
 \epsilon_\infty = \max_k \max_{x\in[0,1]}|H_k^n - H(1,x_k)|,\quad  \epsilon_i = \frac{|H^n - H(1,\cdot)|_{\ell^i}}{|H(1,\cdot)|_{\ell^i}}, \quad i=1,\;2.
\end{equation}
Figure \ref{fig:EMIIMCIP} shows errors in the numerical solutions against mesh size $N^{-1}$.
Grid refinement studies confirm that the second-order convergence in time and space is achieved. The second-order accuracy in the approximation of $\hat{\varepsilon}$ and $\hat{\mu}$ results in the second-order convergence in the numerical solutions.

\section{Conclusion}
We have developed a numerical scheme for one-dimensional hyperbolic equations with variable coefficient.
The method is based on the backward characteristic method and uses the solution and its derivative as unknowns and cubic Hermite interpolation for each computational cell.
The consistency and the conditional stability of the method was presented.
We have proposed a numerical scheme for one-dimensional hyperbolic
equations in a discontinuous media. We have constructed the immersed interface cubic polynomial.
We have extended the method to the one-dimensional Maxwell's equations with variable material properties by approximating with a piecewise constant media.
\section{Acknowledgement}
This research was supported in part by the Air Force Office of Scientific
Research under grant number FA9550-09-1-0226.
\section{Appendix}
{\mbox{ }}
\\
{\bf{ Proof of Lemma \ref{lem:rhosimple}}} \\
Let $p(z)$ be the characteristic polynomial of $G_{\theta,\lambda}$.  $p(z)$ is written as
\begin{align*}
p(z)&= z^2 +2\left( e^{-i \theta }\lambda  \left(1-3 \lambda
+\lambda ^2\right)- \left(1-2 \lambda +\lambda ^3\right)\right)z
+(1-\lambda )^4+e^{-2 i \theta }\lambda ^4-2 e^{-i \theta } \lambda
\left(1-2 \lambda ^2+\lambda ^3\right) \\
& = z^2 + \beta z + \gamma.
\end{align*}
The necessary and sufficient condition for $\rho_{1,\theta,\lambda} \neq \rho_{2,\theta,\lambda}$ is that $\beta^2-4\gamma \neq 0$. We obtain, after some work,
\begin{equation*}
 \beta^2-4\gamma =
4 (-1+\lambda )^2 \lambda ^2 (\cos \theta  -i \sin \theta) [2 \left(5+\lambda -\lambda ^2\right)+\left(-1-2 \lambda +2 \lambda ^2\right) \cos\theta+3 i (-1+2 \lambda ) \sin\theta],
\end{equation*}
and so $\beta^2-4\gamma \neq 0$ is equivalent to
\begin{equation*}
q(\lambda,\theta):=(2 \left(5+\lambda -\lambda ^2\right)+\left(-1-2 \lambda +2 \lambda ^2\right) \cos\theta)^2+ (3 (-1+2 \lambda ) \sin\theta)^2\neq 0.
\end{equation*}
The last term $(-1+2 \lambda ) \sin\theta$ equals to 0 when $\lambda=1/2$ or $\sin\theta = 0$.
But $q(1/2,\theta) = 3(7-\cos\theta)/2 \neq 0$, $q(\lambda,0)=9$ and $q(\lambda,\pi)= 11 + 4\lambda - 4\lambda^2$. Thus we see that $q(\lambda,\theta) \neq 0$ for all $0 < \lambda < 1$.
$\square$ \\
{\bf {Proof of Lemma \ref{lem:rho_less_thanone}}.}\\
We employ the theory of Schur to check whether the roots of the polynomial $p(z)$ reside inside the unit circle.
Let $p^\ast(z)$ and $p_1(z) $ be polynomials defined by $p^\ast(z) := \bar{\gamma}z^2 + \bar{\beta}z + 1$ and $ p_1(z):=\frac{p^\ast(0)p(z) - p(0)p^\ast(z)}{z} = (1 - |\gamma|^2)z+(\beta-\gamma\bar{\beta})$  respectively.
From Theorem 4.3.2 in \cite{Strikwerda-Finidiffschepart:89}, the eigenvalues $\rho_i$ of $p(z)$ satisfy $|\rho_1| < 1$ and $|\rho_2| < 1 $ if and only if $|p(0)| < |p^\ast(0)|$ and the zero of $ p_1(z)$, which we denote by $\eta$, satisfies $|\eta| < 1$.
The first inequality is equivalent to the inequality $1 > |\gamma|$, and the second one $|\eta| < 1$ is equivalent to $|1-|\gamma||^2 > |\bar{\beta} - \beta \bar{\gamma}|^2$.\\
We obtain after some works
\begin{align*}
&1-|\gamma|^2
%&= 2 \kappa  \left(4-12 \kappa +4 \kappa ^2-3 \kappa ^3+\left(2-6 \kappa -4 \kappa ^2
%+4 \kappa ^3\right) \cos\theta-\kappa ^3 \cos2 \theta \right)\nonumber \\
=4 \kappa  \left(2-6 \kappa +2 \kappa ^2- \kappa ^3+\left(1-3 \kappa- 2 \kappa ^2+2 \kappa ^3\right)
\cos\theta-\kappa ^3 \cos^2 \theta \right)\nonumber \\
&=4 \kappa f(\kappa,\theta),\\
&|1-|\gamma||^2 - |\bar{\beta} - \beta \bar{\gamma}|^2 \nonumber\\
&=(2\kappa\sin\tfrac{\theta }{2})^4 \left(3-12\kappa+11 \kappa^2 -2 \kappa ^3+ \kappa ^4
-2 \kappa ^2 \left(1-\kappa +\kappa ^2\right) \cos\theta+2\kappa ^4 \cos^2\theta \right) \nonumber \\
&=(2\kappa\sin\tfrac{\theta }{2})^4g(\kappa,\theta),
\end{align*}
where $\kappa = \lambda(1-\lambda)$. It is straightforward to see that $f(\kappa,\theta)>0$ and
$g(\kappa,\theta)>0$ for all $0<\kappa\le\tfrac{1}{4}$ and $0<\theta<2\pi$.
Indeed,
\begin{equation*}
 \partial_\kappa f(\kappa,\theta) =
 -3 (\cos\theta-1)^2 \kappa^2 + 4(1-\cos\theta) \kappa
-3(2+\cos\theta)\le -\tfrac{14}{3}-3\cos\theta<0.
\end{equation*}
Hence $f$ is monotone decreasing with respect to $\kappa$ for all
$\theta$. Thus, the inequality
\begin{equation*}
  f(\tfrac{1}{4},\theta) = \tfrac{39 + 10\cos\theta -\cos^2\theta}{64} > 0,
\end{equation*}
implies that $f(\kappa,\theta)>0$. Finally
\begin{align*}
  \partial_{\kappa,\kappa} g(\kappa,\theta)& =
  12(2-2\cos\theta +\cos2\theta)\kappa^2 + 12(\cos\theta -1)\kappa - 4\cos\theta + 22 \\
  &\ge \partial_{\kappa,\kappa} g(\tfrac{1-\cos\theta}{2(2-2\cos\theta-\cos2\theta)},\theta)
  =\tfrac{87-96\cos\theta - 49\cos2\theta + 4\cos3\theta}{2(2-2\cos\theta-\cos2\theta)}>0.
\end{align*}
thus, $\partial_{\kappa}g(\kappa,\theta)$ is increasing with respect
to $\kappa$ for all $\theta$, and
\begin{equation*}
  \partial_\kappa g(\tfrac{1}{4},\theta)=\tfrac{-108-12\cos\theta + \cos2\theta}{16}<0.
\end{equation*}
Therefore $g(\kappa,\theta)$ is decreasing with respect to $\kappa$
for all $\theta$, and we have
\begin{equation*}
  g(\kappa,\theta)\ge g(0,\theta) = \tfrac{170-26\cos\theta + \cos2\theta}{256}>0.
\end{equation*}
%we have
%\[
%|1-|\gamma||^2 - |\bar{\beta} - \beta \bar{\gamma}|^2 =
%(2\kappa\sin\tfrac{\theta }{2})^4g(\kappa,\theta) > 0
%\] for $\theta\in(0,2\pi)$. The proof is completed.
$\square$
\bibliographystyle{siam}
%\bibliography{CIP}

\end{document}